\documentclass[10pt,a4paper]{article}

\usepackage[utf8]{inputenc}
\usepackage{hyperref}
\usepackage[sc]{mathpazo}
\usepackage[T1]{fontenc}
\usepackage{lmodern}        
\usepackage[a4paper]{geometry}
\usepackage{authblk}
\usepackage{amsmath, amsfonts, amssymb}
\usepackage{graphicx}
\usepackage{xcolor}
\usepackage{ntheorem}
\usepackage{sectsty}
\usepackage[english]{babel}

\numberwithin{equation}{subsection} 
 
\theoremstyle{break}
\newtheorem{theorem}{Theorem}
\newtheorem{corol}{Corollary}
\newtheorem{lemme}{Lemma}
\newtheorem{prop}{Proposition}
\newtheorem{definit}{Definition}
\newtheorem{rem}{Remark}

\newcommand{\B}{\mathcal{B}}
\newcommand{\U}{\mathcal{U}}
\newcommand{\R}{\mathcal{R}}
\newcommand{\F}{\mathcal{F}}
\newcommand{\tj}{\tilde{\jmath}}
\newcommand{\W}{\mathcal{W}}
\newcommand{\Cun}{\mathcal{C}^{1,1}}

\allsectionsfont{\centering\mdseries\scshape}

\title{A regularity result for fixed points, with applications to linear response}
\date{\today}
\author{Sedro Julien\footnote{Universit\'e Paris-Sud, Campus d'Orsay, B\^atiment 430, Bureau 108. Contact: Julien.Sedro@math.u-psud.fr}}

\begin{document}
\maketitle
\begin{abstract}
	In this paper, we show a series of abstract results on fixed point regularity with respect to a parameter. They are based on a Taylor development taking into account a \emph{loss of regularity} phenomenon, typically occurring for composition operators acting on spaces of functions with finite regularity. We generalize this approach to higher order differentiability, through the notion of an \emph{n-graded family}.
	\\We then give applications to the fixed point of a non linear map, and to linear response in the context of (uniformly) expanding dynamics (theorem \ref{linearresponse} and corollary\ref{Gmeasurediff}), in the spirit of Gouëzel-Liverani. 
\end{abstract}

\section{Introduction}
The aim of this paper is to study the following regularity problem for a fixed point depending on a (multi-dimensional) parameter : 
\\Given three Banach spaces $\mathcal{B}$, $X_0,~X_1$, such that there exists a continuous, linear injection $j_0:X_1\hookrightarrow X_0$, we consider maps $F_i:\mathcal{U}\times A_i\longrightarrow A_i~(i\in\{0,1\})$, where $\mathcal{U}\subset\mathcal{B}$ open, $A_1\subset X_1$ is closed and non-empty and $A_0=j_0(A_1)$. We assume that for every $\phi_1\in A_1$, every $u\in \mathcal{U}$, $j_0\circ F_1(u,\phi_1)=F_0(u,j_0(\phi_1))$and the existence, for every $u\in \mathcal{U}$, of a $\phi_1(u)\in A_1$, such that

\begin{equation}
F_1(u,\phi_1(u)) = \phi_1(u)
\end{equation}

We consider situations where the fixed-point map $\phi_1:u\in \mathcal{U}\mapsto\phi_1(u)\in X_1$ has no particular regularity, yet when one takes the injection $\phi_0=j_0\circ\phi_1$, one gains some regularity.
\medskip

When studying the regularity of fixed point map, the most natural tool at our disposal is the implicit function theorem, formulated in the Banach space setting. However, there are a number of cases where this approach fails, notably when the maps $F_i$ are not continuously differentiable in the classical sense : this is the case when, for example, F is a composition operator. 
\medskip

One can give explicit examples, where some sort of regularity can be recovered, and much can be obtained through elementary methods: for instance, given $\epsilon>0$, $u\in[-\epsilon,\epsilon]$ and $g\in C^1([-1,1]\times[-\epsilon,\epsilon])$ non identically zero, the operator $F(u,\phi)(t)=\frac{1}{2}\phi(\frac{t+u}{2})+g(t,u)$ initially defined on $[-\epsilon,\epsilon]\times C^1([-1,1])$, is a contraction in its second variable~\footnote{i.e for every fixed $u\in[-\epsilon,\epsilon]$, $||F(u,\phi)-F(u,\psi)||_{C^\alpha}\leq k_u||\phi-\psi||_{C^\alpha}$ with $\sup_{u\in[-\epsilon,\epsilon]}k_u<1$} when acting on $C^\alpha([-1,1])$, for every $\alpha\in[0,1]$.
\medskip

However, the fixed point map, $u\in[-\epsilon,\epsilon]\longmapsto\phi_u\in C^\alpha([-1,1])$ is not continuous. Yet, the map  $u\in[-\epsilon,\epsilon]\longmapsto\phi_u\in C^0([-1,1])$ is $\alpha$-Hölder~(see appendix \ref{elementary}) 
\\Another kind of problem arises when one studies the differentiability of that fixed point map: if it is natural to define a partial differential $\partial_uF(u,\phi)=\frac{1}{4}\phi'(\frac{u+.}{2})+\partial_ug(.,u)$ for every $\phi\in C^1([-1,1])$, the corresponding operator $\partial_uF(u,.)$ does not define a partial differential with respect to $u$ for $F(u,.)$ in the classical sense (as it is not a linear map from $\mathbb{R}$ to $C^1([-1,1])$): a phenomenon one can consider as \emph{loss of regularity}.
\medskip

Our main result, theorem \ref{mainresult}, allows one to obtain differentiability results for the kind of fixed points problems outlined in this introduction. The full statement is given in §~\ref{mainresultstatement} and a proof in §~\ref{proofmainresult}; it is based on a type of Taylor development \eqref{NewTaylor}, which can be interpreted as an analogue of the Gouëzel-Liverani spectral perturbation result (\cite[~§8.1]{GL06}). The major improvement here is the possible application to fixed points of non-linear maps.
\\We also discuss a generalization to higher order differentiability in §~\ref{mainresult2section}, by introducing a notion of \emph{graded family} (definition~\ref{gradedfamily}). 
\medskip

In section \ref{nonlinear} we propose an application of our result to a non-linear situation, where the set of parameter lies in an infinite dimensional space: in short, we interpret the perturbation itself as a parameter, and study regularity of the solutions with respect to it. This example is somehow "minimal", in the sense that it is the simplest non trivial, non-linear example we could think of.
\medskip

We then turn to an application to linear response for expanding dynamics, i.e differentiability results and first-order variations for the absolutely continuous, invariant measure for a one parameter family of dynamics. This field has already been thoroughly studied, in various dynamical contexts: uniformly expanding maps on the circle (\cite{Ba14}) or on general Riemann manifold (\cite{Ba16}), intermittent maps of the interval (\cite{BaladiTodd16}), piecewise expanding maps of the interval (\cite{BaladiSmania08,BaladiSmaniaalt10,BaladiSmania12}).

In the hyperbolic case, one can cite Ruelle's work on Axiom A systems (\cite{Ru97,Ru97erratum}, see also the erratum by Jiang \cite{Jiang12}), Gouëzel-Liverani papers on Anosov and Axiom A systems (\cite{GL06,GL08}), and the 2004 paper of Dolgopyat on partially hyperbolic systems (\cite{Dol04}). In a different vein, one can see the paper by Haider and Majda \cite{HairerMajda10}.
\medskip

The "modern" approach to linear response is based on the "weak spectral perturbation" techniques developed in Gouëzel and Liverani papers (see Baladi's monograph \cite{Ba16} and the original papers \cite{GL06,GL08}, see also \cite{G10}). Our method allows to recover similar regularity results, as well as a linear response formula, and one can fruitfully compare our main results theorem~\ref{mainresult} \ref{spectraldataisholder}, \ref{linearresponse}, and corollary \ref{Gmeasurediff} to Gouëzel-Liverani's paper (\cite[§8.1,~(8.3)]{GL06}),  and to Gouëzel's paper \cite[Cor~3.5,~p.21]{G10} (see also theorem 2.36 and 2.38 in Baladi's book \cite{Ba16}). Let us emphasize the differences and similarities:
\begin{itemize}
	\item \textbf{Linear versus non-linear}: If the "weak spectral perturbation" theorem only applies to bounded, linear operators, our theorem \ref{mainresult} can also be applied to non-linear maps (see §~\ref{nonlinear}). However, it is worth noting that when one does apply our theorem to (linear) transfer operators, the "Taylor development" \eqref{NewTaylor} becomes (8.3) in §8.1 of \cite{GL06} (i.e the Taylor expansion assumption in Gouëzel-Liverani paper): this is made precise in the proof of lemma \ref{perturbedtransferop}.
	
	\item \textbf{Parameter dimension}: Our result is naturally formulated for a multi-dimensional (even infinite-dimensional) parameter, whereas Gouëzel-Liverani spectral theorem assume a one-dimensional parameter. Nonetheless, the latter can easily be extended to multi-dimensional parameter. It is not known whether it can be generalized to an infinite-dimensional parameter. We provide an application with an infinite-dimensional parameter in §~\ref{nonlinear}. 
	
	\item \textbf{Uniform Lasota-Yorke versus fixed point continuity}: The proper generalization of the uniform Lasota-Yorke inequalities (assumptions (8.1-8.2), \cite{GL06}) in Gouëzel-Liverani result seems to be the continuity of the fixed point map. A notable difference in our approaches is that the spectral gap assumption is made on the largest Banach spaces, whereas our fixed point map existence and continuity assumption ((i) in theorem \ref{mainresult}) is on the smallest one. Otherwise, the scheme works in the same sense, i.e gain of one derivative when going to the next space.
	
	\item \textbf{Regularity results for the normalized eigenfunction}: \cite[Cor~3.5,~p.21]{G10} studies the regularity of the normalized eigenfunction $\phi_t$ of a transfer operator $(\mathcal{L}_t)_{t\in(-\delta,\delta)}$. It is shown that when the transfer operator acts on $X_0\hookrightarrow X_1$ 2 Banach spaces with a Taylor expansion of the form \eqref{GLTaylorexpansion}, then $\phi_t$ admits itself a Taylor expansion at $t=0$ in $X_0$: $||\phi_t-\phi_0-t.\partial_t{\phi_t}_{\left |t=0\right.}||_{X_0}=\mathcal{O}(|t|^{2-\epsilon})$ with $\epsilon>0$ arbitrarily small, not depending on the spaces $X_0,~X_1$. We obtain a very similar result in theorem \ref{linearresponse}, by applying theorem \ref{mainresult}.
\end{itemize}
\medskip

In order to keep the exposition to a reasonable length, we will not discuss applications of theorem \ref{thetheorem} to higher-order response theory, nor to higher-order differentiation of the spectral data of the transfer operator. To the reader interested by this subject, we recommend \cite[§8.1]{GL06} or \cite{Ru98}.
\medskip

A fair warning to our reader : throughout the text, constants are denoted by the letter C, whose numerical value changes from one occurrence to the next.
\medskip

Recall that if $\Omega\subset\mathbb{R}^n$ is an open subset, $f\in C^0(\Omega)$, $k\in\mathbb{N}$, $\alpha\in(0,1)$, $r=k+\alpha>0$, we say that $f$ is a $C^{r}$ map on $\Omega$ if f is of class $C^k$ on $\Omega$ and its k-th differential (seen as a k multi-linear map) is $\alpha$-Hölder.
We endowed the space of $C^{r}$ maps of $\Omega$ with the norm

\begin{equation}
||f||_{C^{r}}=\max(||f||_{C^{k}}, \underset{x\not=y}{\sup}\frac{||D^kf(x)-D^kf(y)||}{||x-y||^\alpha})
\end{equation}
\medskip
 
The author would like to thank the anonymous referees for their many suggestions and comments, which greatly helped improve both the presentation and mathematical content of the paper.   
\\The author also acknowledges the support of the ESI in Vienna, where the redaction of this paper was started in May 2016.
\\Finally, the author would like to express his warmest thanks to Hans Henrik Rugh, for his constant support, his availability, and many fruitful conversations during the maturation of this work.

\section{Differentiation and graded diagram}
\subsection{Main results}\label{mainresultstatement}

This theorem can be thought of as a complement to the implicit function theorem. Besides the resemblance with \cite[Thm 8.1]{GL06} one can see an analogy with the Nash-Moser scheme \cite{hamilton1982inverse}, with the use of a (finite) scale of spaces. 
\begin{definit}[Scale of Banach spaces]\label{scaledef}
Let $n\geq 1$. A family of Banach spaces $X_{0}\supset X_{1}\supset\dots\supset X_{n}$ is said to be a \emph{scale} if the injective linear maps $j_k:X_{k+1}\rightarrow X_k$ are bounded (i.e $0\leq i\leq j\leq n\Leftrightarrow~||.||_{X_j}\leq||.||_{X_i}$).
\\We will denote a scale by $X_{0}\overset{j_0}{\hookleftarrow}X_{1}\overset{j_1}{\hookleftarrow}\dots\overset{j_{n-1}}{\hookleftarrow}X_n$, or simply by $(X_n,\dots,X_1,X_0)$.
\end{definit}
Note that scales of spaces already appeared in \cite{G10,GL06} and other previous works on spectral stability (\cite{BY93,KL98}).

\begin{theorem}\label{mainresult}
	Let $\mathcal{B}$, $X_{0}, X_{1}$ be Banach spaces such that $X_{0}\overset{j_0}{\hookleftarrow} X_{1}$.
	\\Let $A_{1}\subset X_{1}$ be closed and non-empty, and $A_{0}=j_{0}(A_{1})\subset X_{0}$. 
	\\Let $u_{0}\in\mathcal{B}$, and $\mathcal{U}$ a neighborhood of $u_{0}$ in $\mathcal{B}$. 
	\\Consider continuous maps $F_{i}:\mathcal{U}\times A_{i}\rightarrow A_{i}$, where $i\in\{0,1\}$, with the following property :
	\begin{equation}
	F_{0}(u,j_{0}(\phi_{1}))=j_{0}\circ F_{1}(u,\phi_{1})
	\end{equation}
	for all $u\in \mathcal{U}$, $\phi_{1}\in A_{1}$.
	
	Moreover, we will assume that :
	
	\begin{enumerate}
		
		\item [(i)] For every $u\in \mathcal{U}$, $F_{1}(u,.):A_{1}\rightarrow A_{1}$ admits a fixed point $\phi_{1}(u)\in A_{1}$.
		\\Furthermore, the map $u\in \mathcal{U}\longmapsto\phi_{1}(u)\in X_{1}$ is continuous. 
		\item [(ii)]Let $\phi_0(u)=j_0(\phi_1(u))$. 
		\\For some $(u_{0},\phi_{0}(u_0))=(u_0,\phi_0)\in \mathcal{U}\times j_{0}(A_{1})$, there exists $P_0=P_{u_{0},\phi_{0}}\in L(\mathcal{B},X_{0})$, $Q_0=Q_{u_{0},\phi_{0}}\in L(j_{0}(X_{1}),X_{0})$, such that
		\begin{equation}\label{NewTaylor}
		F_{0}(u_{0}+h,\phi_{0}+z_{0})-F_{0}(u_{0},\phi_{0})=P_0.h+Q_0.z_{0}+(||h||_{\mathcal{B}}+||z_{0}||_{X_{0}})\epsilon(h,z_{1})
		\end{equation}
		where $h\in\mathcal{B}$ satisfies $u_{0}+h\in \mathcal{U}'$, $z_{1}\in A_{1}$, $z_{0}=j_{0}(z_{1})\in A_{0}$, and $\epsilon(h,z_{1})\overset{X_{0}}{\underset{(h,z_{1})\rightarrow(0,0)}{\longrightarrow}} 0$.
		\item [(iii)]$Id-Q_0\in L(j_{0}(X_{1}),X_{0})$ can be extended to a bounded, invertible operator of $X_{0}$ into itself.
	\end{enumerate}
	
	Then the following holds :
	
	\begin{enumerate}
		\item [(i)']Let $\phi_{0}(u)=j_{0}(\phi_{1}(u))$. The map $u\in \mathcal{U}\mapsto \phi_{0}(u)\in X_0$ is differentiable at $u=u_{0}$
		\footnote{i.e there exists a bounded, linear operator $D_u\phi_{0}(u_0):\mathcal{B}\rightarrow X_0$ such that $||\phi_0(u_0+h)-\phi_0(u_0)-D_u\phi_0(u_0).h||_{X_0}\underset{h\rightarrow 0}{\rightarrow}0$ for all $h\in\mathcal{B}$ such that $u_0+h\in \mathcal{U}$.}. 
		\item [(ii)'] Its differential satisfies 
		\begin{equation}
		D_u\phi(u_{0})=(Id-Q_0)^{-1}P_0
		\end{equation}
	\end{enumerate}
\end{theorem}

\begin{rem}
If one were to take $\epsilon(h,z_1)$ in \eqref{NewTaylor} depending only upon $h$, one could recover a condition similar to \cite[§8.1,~(8.3)]{GL06} (see lemma \ref{perturbedtransferop})
\end{rem}	

It can seem artificial to include a statement about continuity of the map $u\in \mathcal{U}\longmapsto \phi_{1}(u)\in X_{1}$ without further explanation. One of the cases where such an assumption can be rigorously justified is when one of the iterates of $F_{1}:\mathcal{U}\times A_{1}\rightarrow A_{1}$, say $F_1^n$ is a \emph{contraction} w.r.t its second variable, a classical result in fixed point theory:

\begin{prop}\label{fixedpointcontinuity}
Let $\mathcal{B},~X$ be Banach spaces, $\U\subset\mathcal{B}$ an open set and $A\subset X$ be closed, non-empty. Let $F:\U\times A\rightarrow A$ be a continuous map, such that there exists $n\in\mathbb{N}$ for which $F^n$ is a contraction with respect to its second variable.
\\Then for every $u\in\U$, $F(u,.)$ admits a unique fixed point $\phi_u\in A$, and furthermore the map $u\in\U\mapsto\phi_u\in X$ is continuous. 
\end{prop}	

\paragraph{Proof of proposition \ref{fixedpointcontinuity}:}
We can apply the Banach contraction principle to $F^n:\mathcal{U}\times A\rightarrow A$, and thus obtain the existence of a fixed point $\phi(u)\in A$ for every $u\in\mathcal{U}$. We also have :
	
	\begin{align}
	||\phi(u)-\phi(u_{0})||_{X}&=||F^n(u,\phi(u))-F^n(u_{0},\phi(u_{0}))||_{X}\\
	&=||F^n(u,\phi(u))-F^n(u_{0},\phi(u))+F^n(u_{0},\phi(u))-F^n(u_{0},\phi(u_{0}))||_{X}\\
	&\leq C||\phi(u)-\phi(u_{0})||_{X}+ ||F^n(u,\phi(u))-F^n(u_{0},\phi(u))||_{X}
	\end{align}
	
	with $C<1$, so that :
	
	\begin{equation}
	||\phi(u)-\phi(u_0)||_{X}\leq \dfrac{1}{1-C}||F^n(u,\phi(u))-F^n(u_0,\phi(u))||_{X}
	\end{equation}
	
	We can now conclude with the continuity of $F:\mathcal{U}\times A\rightarrow A$.
\medskip
 
Remark that if we were to demand a stronger condition on the regularity of $F$ with respect to $u$, say Hölder-continuity or Lipschitz continuity, the fixed point map $u\in \mathcal{U}\rightarrow\phi(u)\in X$ would mirror that condition.
\medskip

In section \ref{mainapp}, we illustrate the abstract theorem \ref{mainresult} by applying it to a positive, linear transfer operator $\mathcal{L}_u$, associated with a family $(T_u)_{u\in \mathcal{U}}$ of $C^{1+\alpha}$ expanding maps on a Riemann manifold X, acting on $C^{1+\alpha}(X)$, and who admits an isolated, simple eigenvalue $\lambda_u$ of maximal modulus. 
It requires to work with the nonlinear map $F:\mathcal{U}\times C^{1+\alpha}(X)$, defined for $u\in \mathcal{U}$ a neighborhood of $u_0\in\mathcal{B}$ and $\phi\not\in\ker\mathcal{L}_{u}^{*}\ell_{u_{0}}$, by
\begin{equation}\label{normalisation}
F(u,\phi)= \dfrac{\mathcal{L}_u\phi}{\int\mathcal{L}_u\phi d\ell_{u_0}}
\end{equation}
where $\ell_u$ (resp. $\phi_u$) is the left (resp. right) eigenvector of $\mathcal{L}_u$, chosen so that $\int\mathcal{L}_u\phi_ud\ell_{u}=\lambda_u$. For $u\in\U$, we chose $\phi_u$ so that $\langle\ell_{u_0},\phi_u\rangle=1$ (this will prove useful in §~\ref{holderspectraldatasection}).

This (nonlinear) renormalization originates from cone contraction theory, and has been used e.g in \cite{Rugh08,Rugh10}. Satisfying  assumption (iii) in theorem \ref{mainresult} is the main reason why one is lead to introduce $\eqref{normalisation}$: indeed, working with the naive guess $\lambda_u^{-1}\mathcal{L}_u$ (for which $\phi_u$ is an obvious fixed point) cannot give a bounded and invertible second partial differential, by definition of an eigenvalue... 
\\It is also worth noting that the normalized maps F satisfy condition (i) in theorem \ref{mainresult} thanks to proposition \ref{fixedpointcontinuity}. More precisely, we are able to establish the following:

\begin{theorem}\label{spectraldataisholder}
For every $0\leq\beta<\alpha$, $u\in \mathcal{U}$, one has
	\begin{itemize}
		\item F(u,.) acts continuously (and even analytically) on $C_+^{1+\alpha}(X)^*:=\{f\in C^{1+\alpha}(X), f\geq 0~and~f\not=0\}$ 
		\item Consider $F(u,.):C_+^{1+\alpha}(X)^*\longmapsto C_+^{1+\beta}(X)^*$. Then $u\in \mathcal{U}\mapsto F(u,.)$ is Hölder continuous, with exponent $\gamma:=\alpha-\beta$.
		\item $F(u,.)$ admits a unique fixed point $\phi(u)\in C_+^{1+\alpha}(X)^*$, and $u\in \mathcal{U}\longmapsto\phi(u)\in C^{1+\beta}(X)$ is $\gamma$-Hölder.  
	\end{itemize}
\end{theorem}
We establish this result in §\ref{holderspectraldatasection}. It also establishes the first assumption of theorem \ref{mainresult}, and is therefore instrumental in proving the following:
 
\begin{theorem}\label{linearresponse}
	Let $0\leq\beta<\alpha<1$, $u_{0}\in\mathcal{B}$, $\mathcal{U}$ a neighborhood of $u_{0}$, $(T_{u})_{u\in \mathcal{U}}$ be a family of $C^{1+\alpha}$, expanding maps of a Riemann manifold X. For each $u\in \mathcal{U}$, let $\mathcal{L}_{u}$ be a weighted transfer operator on $C^{1+\alpha}(X)$, associated with $T_{u}$, defined by \eqref{transferopdef}. 
	\\Let $\lambda_{u}>0$ be its dominating eigenvalue, $\phi(u)\in C^{1+\alpha}(X)$, $\ell_{u}\in (C^{1+\alpha}(X))^{*}$ be the associated eigenvectors of $\mathcal{L}_{u}$ and $\mathcal{L}_{u}^{*}$ respectively. We denote by $\Pi_u$ the associated spectral projector, and let $R_u=\mathcal{L}_u-\lambda_u\Pi_u$ (cf appendix \ref{expmapspectrum}).
	\\Then the following holds true:
	\begin{itemize}
		\item The map $u\in \mathcal{U}\longmapsto \phi(u)\in C^{\beta}(X)$ is differentiable.
		\item We have the following linear response formula for the derivative with respect to u at $u=u_0$:
		\begin{equation}\label{linearresponseformula}
		D_{u}\phi(u_{0})=\dfrac{1}{\lambda_{u_0}}(Id-\lambda_{u_{0}}^{-1}R_{u_{0}})^{-1}(Id-\Pi_{u_0})\partial_{u}\mathcal{L}_{u}|_{u=u_{0}}
		\end{equation}
		
	\end{itemize}
\end{theorem}
We establish this result in §\ref{diffspectraldatasection}, by applying theorem \ref{mainresult} to $F$ acting on the scale $(C^{1+\beta}(X),C^{\beta}(X))$ for any $0<\beta<\alpha$. We show that F satisfies to a Taylor expansion of the form \eqref{NewTaylor}, with (see formulas \ref{partial1},~\ref{partial2}~) 
$$P_0=\frac{1}{\lambda_{u_0}}(Id-\Pi_{u_0})\partial_{u}\mathcal{L}_u|_{u=u_0}\phi_0$$ $$Q_0=\frac{1}{\lambda_{u_0}}\mathcal{L}_{u_0}-\Pi_{u_0}$$ 

\subsection{Taking the first derivative : a proof of theorem \ref{mainresult}}\label{proofmainresult}

Thanks to assumption $(ii)$, we can estimate the difference $z_{0}(h)=\phi_{0}(u_{0}+h)-\phi_{0}(u_{0})$ for $h\in\mathcal{B}$, $u_{0}+h\in \mathcal{U}$. 

\begin{align*}
\phi_{0}(u_{0}+h)-\phi_{0}(u_{0})&=F_{0}(u_{0}+h,\phi_{0}(u_{0}+h))-F_{0}(u_{0},\phi_{0})\\
&=F_{0}(u_{0}+h,\phi_{0}(u_{0})+z_{0})-F_{0}(u_{0},\phi_{0})\\
&=P_0.h+Q_0.z_{0}(h)+(||h||_{\mathcal{B}}+||z_{0}(h)||_{X_{0}})\epsilon(h,z_{1})
\end{align*}

thus, by $(iii)$:

\begin{equation}\label{Increment}
z_{0}(h)=(Id-Q_0)^{-1}P_0.h+(Id-Q_0)^{-1}(||h||_{\mathcal{B}}+||z_{0}(h)||_{X_{0}})\epsilon(h,z_{1})
\end{equation}

Now, remark that :

\begin{itemize}

\item By continuity of $u\in \mathcal{U}\rightarrow\phi_{1}(u)\in X_{1}$ (which is assumption $(i)$), we have $\underset{h\rightarrow 0}{\lim}~z_{1}(h)=0$ in $X_1$, so that 
$\epsilon(h,z_{1}(h))=\epsilon(h)\rightarrow 0$ in $X_0$ as $h\rightarrow 0$ in $\mathcal{B}$.

\item $(Id-Q_0)^{-1}\epsilon(h,z_{1})||h||_{\mathcal{B}}=o(h)$ in $X_0$ as $h\rightarrow 0$ in $\mathcal{B}$

\item For $h$ small enough in $\mathcal{B}$-norm,
\begin{equation}
||(Id-Q_0)^{-1}||.||\epsilon(h)||_{X_{0}}\leq \dfrac{1}{2}
\end{equation}

\end{itemize}

Thus, taking the $X_{0}$-norm in \eqref{Increment} and choosing $h$ small enough in $\mathcal{B}$-norm, we obtain :

\begin{align}
\notag||z_{0}(h)||_{X_{0}}&\leq||(Id-Q_0)^{-1}P_0.h||_{X_0} +||(Id-Q_0)^{-1}\epsilon(h,z_{1})||_{X_0}||h||_{\mathcal{B}}+\dfrac{1}{2}||z_{0}(h)||_{X_{0}}\\
\dfrac{1}{2}||z_{0}(h)||_{X_{0}}&\leq||(Id-Q_0)^{-1}P_0.h||_{X_0} +||(Id-Q_0)^{-1}\epsilon(h,z_{1})||_{X_0}||h||_{\mathcal{B}}
\end{align}

and thus :

\begin{equation}\label{Maj}
z_{0}(h)=\mathcal{O}(h)
\end{equation}

Following \eqref{Maj}, the second term of the sum in the right hand term of \eqref{Increment} becomes :

\begin{equation}
(Id-Q_0)^{-1}(||h||_{\mathcal{B}}+\mathcal{O}(h))\epsilon(h)=o(h)
\end{equation}

Finally, in the $X_0$-topology,

\begin{equation}
z_{0}(h)=(Id-Q_0)^{-1}P_0.h+o(h)
\end{equation}

and thus $u\in \mathcal{U}\rightarrow \phi_{0}(u)\in X_{0}$ is differentiable at $u=u_{0}$ and

\begin{equation}
D_u\phi_0(u_{0})=(Id-Q_0)^{-1}P_0
\end{equation}

\subsection{Higher differentiability and graded diagram}\label{mainresult2section}

In order to differentiate the fixed point map we have to consider an argument coming from a smaller, more "regular" space. More precisely, we showed that if there is, for every $u\in \mathcal{U}$ a $\phi_{1}(u)\in X_{1}$ such that $F(u,\phi_{1}(u))=\phi_{1}(u)$, then $u\longmapsto\phi_{0}(u)=j_{0}(\phi_{1}(u))\in j_{0}(X_{1})\subset X_{0}$ is differentiable.
\bigskip

We aim to iterate this approach to differentiate further the fixed point map with respect to the parameter. In order to do so, we define a notion of an \emph{n-graded family} as such :
	\begin{definit}[Graded family]\label{gradedfamily}
	Let $n\geq 1$ be an integer, and consider a Banach space $\mathcal{B}$, a scale $X_{0}\overset{j_0}{\hookleftarrow}X_{1}\overset{j_1}{\hookleftarrow}\dots\overset{j_{n-1}}{\hookleftarrow}X_n$, $\mathcal{U}\subset\mathcal{B}$ an open subset, $A_n\subset X_n$ a closed, non-empty subset. 
	\\For $0\leq k<l<n$, we denote by $j_{k,l}$ the bounded linear map $j_{k}\circ j_{k+1}\circ..\circ j_{l}:X_{l+1}\rightarrow X_{k}$, and by $\tj_k=j_k\circ..\circ j_{n-1}:X_n\rightarrow X_k$.
	\medskip
	
	Define, for $i\in\{0,..,n-1\}$, $A_i=j_{i,n-1}(A_n)$, and continuous maps $F_{i}:\mathcal{U}\times A_{i}\rightarrow A_{i}$, $i\in\{0,..,n\}$ such that :
	\begin{itemize}
		\item[(i)] For every $u\in \mathcal{U}$, $\phi_{i}\in A_{i}$, $j_{i}(F_{i+1}(u,\phi_{i+1}))=F_{i}(u,j_{i}(\phi_{i+1}))$
		\item[(ii)] There exists $(u,\phi_n)\in\U\times X_n$, such that for every $h\in\B$ such that $u+h\in\U$, every $z_n\in X_n$, such that $\phi_n+z_n\in A_n$, for every $1\leq k\leq n$, $F_{n-k}$ satisfies
		\begin{equation}\label{NewTaylor2}
		F_{n-k}(u+h,\phi_{n-k}+z_{n-k})-F_{n-k}(u,\phi_{n-k})=\sum_{\ell=1}^{k}\sum_{(i,j)\atop i+j=\ell}Q^{(i,j)}(u,\phi_{\ell})[h,z_{\ell-1}]+\R_n(h,z_n)
		\end{equation}
		where for every pair $(i,j)$ so that $i+j=\ell$, 
		\begin{itemize}
		\item[$\bullet$] $Q^{(i,j)}(u,\phi_\ell)\in L(\B^i\times X_{\ell-1}^j,X_{n-k})$ is a $\ell$-linear map
		\item [$\bullet$] $\R_n\in C^0(\B\times X_n,X_{n-k})$ is such that $||\R_n(h,z_n)||_{X_{n-k}}=o(||h||_{\B}^{k},~||z_{n-1}||_{X_{n-1}}^{k})$.
		\end{itemize}
	\end{itemize}
	
	We call a family of maps $(F_{i})_{i\in\{0,..,n\}}$ acting on $\mathcal{B}$, $X_{0}\overset{j_0}{\hookleftarrow}X_{1}\overset{j_1}{\hookleftarrow}\dots\overset{j_{n-1}}{\hookleftarrow}X_n$ and satisfying (i)-(ii), an \textbf{n-graded family}.
\end{definit}

\begin{lemme}\label{centralassumptionlemma}
	Let $(F_{i})_{i\in\{0,..,n\}}$ be an n-graded family. 
	\\Then for every $1\leq k\leq n$, whenever $\phi_n\in C^{0}(\mathcal{U},A_{n}),\dots,\phi_{n-k+1}=\tj_{n-k+1}(\phi_n)\in C^{k-1}(\mathcal{U},A_{n-k+1})$, 
	the map $u\in \mathcal{U}\mapsto \tj_{n-k}\circ\F_{n}(u,\phi_{n}(u))\in X_{n-k}$ is $k-1$ times differentiable, and has the following property: there exists $R^{(k-1)}(u)\in L^{k-1}(\mathcal{B},X_{n-k})$ such that for every $u\in \mathcal{U}$
	\begin{equation}\label{centralassumption}
	D_{u}^{k-1}[\tj_{n-k}\circ\F_{n}](u,\phi_{n}(u))= R^{(k-1)}(u)+Q^{(0,1)}_{u,\phi_{n-k+1}(u)}. D^{k-1}_u \phi_{n-k}(u)
	\end{equation}
	
	where 
	\begin{enumerate}
		\item $u\in \mathcal{U}\longmapsto R^{(k-1)}(u)\in L^{k-1}(\mathcal{B},X_{k})$ is differentiable.
		\item $u\in \mathcal{U}\longmapsto Q^{(0,1)}(u,\phi_{n-k+1}(u))\in L(X_{n-k})$ is differentiable.
	\end{enumerate}
\end{lemme}
\paragraph{Proof of lemma \ref{centralassumptionlemma}:}
From \eqref{NewTaylor2}, one can write:
\begin{align}\label{taylorestimate1}
\notag\tj_{n-k}&\circ[\F_n(u+h,\phi_n(u+h))-\F_n(u,\phi_n(u))]=\F_{n-k}(u+h,\phi_{n-k}(u+h))-\F_{n-k}(u,\phi_{n-k}(u))\\
&=\sum_{\ell=1}^{k}\sum_{(i,j)\atop i+j=k-\ell+1}Q^{(i,j)}(u,\phi_{n-\ell+1})[h,\phi_{n-\ell}(u+h)-\phi_{n-\ell}(u)]+\R_n(h,\phi_{n}(u+h)-\phi_{n}(u))
\end{align}

From our assumptions, for every $\ell\in\{1,\dots,k-1\}$, $\phi_{n-\ell}$ is $\ell$ times differentiable on $\U$, so that one can write for every $u\in\U$ and $h\in\B$ such that $u+h\in\U$,
\begin{equation}\label{taylorestimate2}
\phi_{n-\ell}(u+h)-\phi_{n-\ell}(u)=D_u\phi_{n-\ell}(u).h+\dots+D_u^{\ell}\phi_{n-\ell}(u)+o(||h||^{\ell})
\end{equation}
the term in $o(||h||^{\ell})$ being understood in $X_{n-\ell}$. 
\\The same Taylor development (at order $k-1$) holds for $\phi_{n-k}=j_{n-k}(\phi_{n-k+1})$. 
Injecting \eqref{taylorestimate2} in \eqref{taylorestimate1} establishes first that $\tj_{n-k}\circ\F_n(.,\phi_n(.))$ is (k-1) times differentiable.
\medskip

Secondly, from the variety of terms it yields, we only choose the terms that are $k-1$ linear in $h$ : this gives us the k-1 differential with respect to u, written as \eqref{centralassumption}, along with the following explicit expression for $R^{(k-1)}$:
\begin{align}
\notag R^{(k-1)}(u)=\sum_{\ell=2}^{k-2}&\sum_{(i,j),i\not=0\atop i+j=\ell}\sum_{r_1,\dots,r_j\leq\ell+1\atop i+r_1+\dots+r_j=k-1}Q^{(i,j)}(u,\phi_{n-\ell}(u))\left[h,D_u^{r_1}\phi_{n-\ell-1}(u),\dots,D_u^{r_j}\phi_{n-\ell-1}(u)\right]\\
&+\sum_{(i,j)\atop i+j=k-1}Q^{(i,j)}(u,\phi_{n-1}(u))\left[h,D_u\phi_{n-2}\right]
\end{align}
From there, it easy to check differentiability with respect to $u$, as the previous expression only involves terms of indices $n-\ell$ with at most $\ell-1$ derivatives.

\begin{theorem}\label{thetheorem}
	Let $(F_{i})_{i\in\{0,..,n\}}$ be a n-graded family. Let $u\in \mathcal{U}$. We make the following assumptions :
	\begin{itemize}
		\item For every $u\in \mathcal{U}$, $F_{n}(u,.)$ admits a fixed point $\phi_{n}(u)$. Furthermore, we assume that the map $u\in \mathcal{U}\longmapsto \phi_{n}(u)\in X_{n}$ is continuous.
		\item For every $0\leq k\leq n-1$, $Id-Q^{(0,1)}_{u,\phi_{k}(u)}$ is an invertible, bounded operator of $X_{k}$.
	\end{itemize}
	Then for every $1\leq k\leq n$ the fixed point map $u\in \mathcal{U}\longmapsto\phi_{n-k}(u)=\tj_{n-k}(\phi_{n}(u))\in X_{n-k}$ is k times differentiable, and one has the following formula for its differential:
	\begin{equation}\label{kdiffformula}
	D_u^k\phi_{n-k}(u)=(Id-Q^{(0,1)}(u,\phi_{n-k+1}(u)))^{-1}R^{(k)}(u)
	\end{equation}
	Furthermore, when $u\in\U\mapsto(Q^{(0,1)}(u,\phi_{n-k+1}(u)),R^{(k)}(u))$ is continuous, then so is $u\in\U\mapsto D_u^k\phi_{n-k}(u)$, i.e the fixed point map $\phi_{n-k}$ is $C^k$.
\end{theorem}
\paragraph{Proof of theorem \ref{thetheorem}:}
The continuity statement is obvious. We present a proof by finite and descending induction.
\begin{itemize}
	\item For $k=1$, the differentiability of $u\in \mathcal{U}\longmapsto\phi_{n-1}(u)$ at $u=u_0$ is simply theorem \ref{mainresult}.
	\item For $k=2$ 
	\\For every $h\in \mathcal{B}$, $u_{0}+h\in \mathcal{U}$, one has, thanks to the case $k=n-1$ and assumption \eqref{centralassumption}:
	\begin{align}
	D_{u}(j_{n-2,n-1}\circ F_{n})&(u,\phi_{n}(u)).h=D_{u}F_{n-2}(u,\phi_{n-2}(u)).h\\
	&=Q^{(1,0)}(u,\phi_{n-1}(u)).h+Q^{(0,1)}(u,\phi_{n-1}(u))D_{u}\phi_{n-2}(u).h
	\end{align} 
	Note that $\phi_{n-2}(u)$ is, for every $u\in \mathcal{U}$, a fixed point of $F_{n-2}(u,.)$, so that
	\begin{equation}
	D_{u}F_{n-2}(u,\phi_{n-2}(u).h=D_{u}\phi_{n-2}(u).h
	\end{equation}
	
	This last equality yields,
	\begin{align}
	&D_{u}\phi_{n-2}(u).h=Q^{(1,0)}(u,\phi_{n-1}(u)).h+Q^{(0,1)}(u,\phi_{n-1}(u))D_{u}\phi_{n-2}(u).h\\
	&(Id-Q^{(0,1)})(u,\phi_{n-1}(u))D_{u}\phi_{n-2}(u).h=Q^{(1,0)}(u,\phi_{n-1}(u)).h\\
	&D_{u}\phi_{n-2}(u_{0}).h=(Id-Q^{(0,1)}(u,\phi_{n-1}(u)))^{-1}Q^{(1,0)}(u,\phi_{n-1}(u)).h
	\end{align}
	By \eqref{centralassumption} in definition \ref{centralassumptionlemma}, 
	$$
	\left\{
	\begin{aligned}
	u\in \mathcal{U}&\mapsto Q^{(1,0)}(u,\phi_{n-1}(u))\\
	u\in \mathcal{U}&\mapsto Q^{(0,1)}(u,\phi_{n-1}(u))
	\end{aligned}
	\right.
	$$
	are differentiable at $u=u_{0}$, between the Banach spaces $\mathcal{B}$ and $L(\mathcal{B},X_{n-2})$ (respectively $L(X_{n-2})$).
	\\By the previous equality, $u\in \mathcal{U}\mapsto D_{u}\phi_{n-2}(u)$ is differentiable at $u=u_{0}$, i.e $u\in \mathcal{U}\mapsto \phi_{n-2}(u)\in X_{n-2}$ is a twice differentiable map at $u=u_{0}$.
	\item Let $3\leq k\leq n$. 
	\\Assume the property :
	\begin{center}
		$u\in\mathcal{U}\longmapsto\phi_{n-k+1}(u)=\tj_{n-k+1}(\phi_{n}(u))\in X_{n-k+1}$ is a k-1 times differentiable map.
	\end{center}
	Then, by lemma \ref{centralassumptionlemma}, \eqref{centralassumption} one can write, for the k-1 differential of $u\mapsto \tj_{n-k}\circ F_{n}(u,\phi_{n}(u))$
	\begin{align}
	\notag D_{u}^{k-1}\phi_{n-k}(u)&=D_{u}^{k-1}\tj_{n-k}\circ F_{n}(u,\phi_{n}(u))\\
	&=R^{(k-1)}(u)+Q^{(0,1)}(u,\phi_{n-k+1}(u))D_{u}^{k-1}\phi_{n-k}(u)
	\end{align}
	
	
	Thus, we obtain by the invertibility assumption,
	\begin{equation}
	D_{u}^{k-1}\phi_{n-k}(u)=(Id-Q^{(0,1)}(u,\phi_{n-k+1}(u)))^{-1}R^{(k-1)}(u)
	\end{equation}
	By virtue of lemma \ref{centralassumptionlemma}, one obtains the differentiability of $u\in \mathcal{U}\mapsto D_{u}^{k-1}\phi_{n-k}(u)$, and therefore, that the map $u\in\mathcal{U}\mapsto\phi_{n-k}(u)\in X_{n-k}$ is k times differentiable, with the announced formula.
	
\end{itemize}

\section{A non linear application}\label{nonlinear}

In this section we give an application of theorem \ref{mainresult} to the study of a fixed point of a non linear map. Note also that the parameters lie in an infinite dimensional space.
\bigskip

Consider the interval $I=[-1,1]$, and let $\Cun(I)$ be the set of $C^1$ map on $I$ with Lipschitz derivative, endowed with the norm $||f||_{1,1}=\max(||f||_{C^1},\underset{x,y\in I\atop x\not=y}{\sup}\dfrac{f'(x)-f'(y)}{x-y})$, which makes it a Banach space. Define the map $F:\Cun(I)\times\Cun(I)\rightarrow\Cun(I)$ by
\begin{equation}\label{nonlinearex}
F(u,\phi)=\frac{1}{2}\phi\circ\phi+u
\end{equation}

We will show the following:
\begin{theorem}\label{nonlinearapp}
	Let $I,~\Cun(I)$, and $F:\Cun(I)\times\Cun(I)\rightarrow\Cun(I)$ be as above. One has:
	\begin{itemize}
		\item[(i)] Let $\U=B_{\Cun}(0,r')$ be an open ball in $\Cun(I)$. There is $r,r'\in(0,1)$, such that for every $u\in\U$, $B_{\Cun}(0,r)$,  $F(u,.)$ is a contraction of $B_{\Cun}(0,r)$ in the $C^1$ topology: therefore it admits a fixed point $\varphi_u\in B_{\Cun}(0,r)$, and furthermore the map $u\in\U\mapsto\varphi_u\in C^1(I)$ is continuous.
		\item [(ii)]$F$ acting on the scale $(C^1(I),C^0(I))$ satisfies a development of the form \eqref{NewTaylor}. Therefore the map $u\in\U\mapsto \varphi_u\in C^{0}(I)$ is differentiable.
	\end{itemize}
\end{theorem}
\paragraph{Proof of theorem \ref{nonlinearapp}:}
\begin{itemize}
\item[(i)] It is a straightforward computation: for every $u\in \U$, one has 
	\begin{align*}
	||F(u,\phi)||_\infty&\leq\frac{||\phi||_{\infty}}{2}+||u||_{\infty}\\
	||D_tF(u,\phi)||_\infty&\leq\frac{||\phi'||_{\infty}^2}{2}+||u'||_{\infty}\\
	||D_tF(u,\phi)||_{Lip}&\leq\frac{||\phi'||_\infty.||\phi'||_{Lip}}{2}(1+||\phi'||_\infty)+||u'||_{Lip}
	\end{align*}
Therefore we should choose $r,r'$ such that $\frac{r}{2}+r'\leq r$, $\frac{r^2}{2}+r'\leq r$ and $\frac{r^2}{2}(1+r)+r'\leq r$. This conditions, which admits obvious solutions, insure us that $F(u,.)$ preserves $B_{\Cun}(0,r)$. From now on, we fix $r,r'$ so that those conditions are satisfied.
\medskip

We now show that $||F(u,\phi)-F(u,\psi)||_{C^1}\leq k||\phi-\psi||_{C^1}$, when $\phi,\psi\in B_{\Cun}(0,r)$. It is noteworthy that here, $k$ is independent of $u$. One has:

\begin{align*}
||F(u,\phi)-F(u,\psi)||_\infty&\leq\frac{1}{2}(1+||\phi'||_{\infty})||\phi-\psi||_{C^1}\\
||D_tF(u,\phi)-D_tF(u,\psi)||_{\infty}&\leq\frac{1}{2}(||\psi'||_\infty+|\phi'|_{Lip}||\phi'||_\infty+||\phi'||_\infty)||\phi-\psi||_{C^1}
\end{align*}

so that one need to impose the following conditions on $r$: $\frac{1+r}{2}<1$, $\frac{2r+r^2}{2}<1$. 
\\Not only do these conditions clearly have solutions, they are also compatible with the conditions imposed on $r$ in (i). From now on, we assume that $r,r'$ satisfy both sets of conditions. 
\medskip

Thus, for every $u\in B_{\Cun}(0,r')$, $F(u,.):B_{\Cun}(0,r)\rightarrow B_{\Cun}(0,r)$ is a contraction in the $C^1$ topology. Hence it admits a fixed point $\varphi_u\in B_{\Cun}(0,r)$, and the map $u\in\U\mapsto\varphi_u\in C^{1}(I)$ is continuous (and even Lipschitz) by proposition \ref{fixedpointcontinuity}.

\item One can write, for $u,h\in\Cun(I)$ such that $u,~u+h\in\U$ and $\phi,z\in C^{1}(I)$,
\begin{align*}
F(u+h,\phi+z)-F(u,\phi)=h+\frac{1}{2}[\phi'\circ\phi.z+z\circ\phi]+(z'\circ\phi).z+||z||_{\infty}\epsilon_0(z)
\end{align*}
where $||\epsilon_0(z)||_{\infty}\longrightarrow 0$ as $||z||_{\infty}\longrightarrow 0$. From there it is clear that with:
\begin{align*}
P_{u,\phi}.h&=h\\
Q_{u,\phi}.z&=\frac{1}{2}[\phi'\circ\phi.z+z\circ\phi]\\
\epsilon(h,z_1)&=(z'\circ\phi).z+||z||_{\infty}\epsilon_0(z)=(z'_1\circ\phi).z_0+||z_0||_{\infty}\epsilon_0(z_0)
\end{align*}
F satisfies a development of the form \eqref{NewTaylor}.
\medskip

To conclude, we need to establish the invertibility (and boundedness of the inverse) of $Q_{u,\phi}=Q_{\phi}$ on $C^0(I)$. 
\\It is easy to see that for every $\phi\in B_{\Cun}(0,r)$, $||Q_\phi.z||_{\infty}\leq\frac{1}{2}(1+r)||z||_\infty$, so that $||Q_\phi||_{C^0}<1$ whenever $r<1$ (which is insured by the sets of conditions in (i),~(ii)). Therefore its Neumann series converges in $C^0(I)$, and $Id-Q_\phi$ has a bounded inverse in $C^0(I)$ for every $\phi\in B_{\Cun}(0,r)$.
\end{itemize} 

\section{Application to linear response for expanding maps}\label{mainapp}

As a second application of our main result theorem \ref{mainresult}, we study the \emph{linear response} problem, in the context of smooth uniformly expanding maps.
\\More precisely, our strategy is the following: 
\begin{itemize}
	\item We first show regularity results (Hölder and Lipschitz continuity, differentiability in the sense of \eqref{NewTaylor}) for the transfer operator $\mathcal{L}_u$ acting on Hölder spaces, with respect to $u$: see lemma \ref{perturbedtransferop}
	\item We then establish theorem \ref{spectraldataisholder} by a direct argument (see §~\ref{holderspectraldatasection}).
	\item We finally prove theorem \ref{linearresponse} by applying theorem \ref{mainresult} to the map $F$ defined by \eqref{normalisation}, acting on the scale $(C^{1+\beta}(X),C^{\beta}(X))$ (see §~\ref{diffspectraldatasection}). 
\end{itemize}

\subsection{Perturbations of the transfer operator}

Let $d\geq 1$, $\epsilon>0$, $\mathcal{U}=(-\epsilon,\epsilon)^{d}$, $0<\alpha<1$ and $(T_{u})_{u\in \mathcal{U}}\in C^{1+\alpha}(X)$ be a $C^{1+\alpha}$ family of $C^{1+\alpha}$ expanding maps. For example, $T_u$ can be a $C^{1+\alpha}$ perturbation of an original expanding map $T_0$: by (iii) in proposition \ref{propofexpmaps}, $T_u$ is also expanding for $u\in \mathcal{U}$ small enough. 
\\Let $g:\U\times X\rightarrow\mathbb{R}$ be a $C^{1+\alpha}$ map. For every $u\in \mathcal{U}$, define the associated transfer operators (e.g, on $L^{\infty}(X)$) by

\begin{equation}\label{transferopdef}
\mathcal{L}_{u}\phi(x)=\sum_{y,T_{u}y=x}g(u,y)\phi(y)
\end{equation}

Recall that the spectral features of interest appears when the transfer operator acts on Hölder spaces (cf appendix~\ref{expmapspectrum}).
In the next proposition, we study the regularity of $\mathcal{L}_u$ with respect to the parameter $u$.

\begin{lemme}[Regularity of the perturbed transfer operator]\label{perturbedtransferop}
Let $0\leq\beta<\alpha<1$, and $\gamma:=\alpha-\beta$. Let $X,\mathcal{U}$ and $g,T_u,\mathcal{L}_u$ be as above.
\begin{itemize}
\item $u\in\mathcal{U}\longmapsto\mathcal{L}_{u}\in L(C^{1+\alpha}(X),C^{1+\beta}(X))$ is $\gamma$-Hölder.
\\In particular, it is a continuous map.
\item For every $h\in\mathcal{B}$ such that $u_{0}+h\in \mathcal{U}$, every $0\leq\beta\leq\alpha$, we can define a bounded operator $\partial_{u}\mathcal{L}(u_0,.).h:C^{1+\beta}(X)\rightarrow C^{\beta}(X)$, such that for every $\phi\in C^{1+\beta}(X)$, 
\begin{equation}\label{GLTaylorexpansion}
\mathcal{L}(u_{0}+h,\phi)-\mathcal{L}(u_{0},\phi)-\partial_{u}\mathcal{L}(u_{0},\phi).h=||h||_{\mathcal{B}}\epsilon(h)
\end{equation}
with $\epsilon(h)\underset{h\rightarrow 0}{\longrightarrow}0$ in $C^{\beta}(X)$
\\Furthermore, $\mathcal{L}$ satisfies \eqref{NewTaylor} in theorem \ref{mainresult}, with the scale $(C^{1+\beta}(X),C^{\beta}(X))$.
\end{itemize}
\end{lemme}

\textbf{Proof:}~By a standard argument (see \cite{Ru89,GuLa03}), one can construct a family of open sets covering X, small enough to be identified with open sets in $\mathbb{R}^{\dim(X)}$, and such that on each of this open sets, the transfer operator is a (finite) sum of operators of the form $\W_u\phi:=(g_u.\phi)\circ\psi_u$, with $\phi\in C^{1+\alpha}(W)$, $\psi\in C^{1+\alpha}(\U\times V,W)$ is a contraction in its second variable (and a local inverse branch of $T_u$), $g\in C^{1+\alpha}(\U\times W)$ with compact support, and $V,~W$ open sets in $\mathbb{R}^{\dim(X)}$.
\\We will apply the results of appendix \ref{compestimates} to the operators $\W_u$. 

\bigskip

For the first item, one needs to estimate, for $\phi\in C^{1+\alpha}(W)$, $||(\W_u-\W_v)\phi||_{C^{1+\beta}}=\max(||(\W_u-\W_v)\phi||_{C^{1}},||D_{x}(\W_u-\W_v)\phi||_{C^\beta})$.

Assume first that the weight $g$ is independent of the parameter. Then by lemma \ref{compisHölder2lemma},~\eqref{compisHölder2}, 

\begin{align}
||(\W_u-\W_v)\phi||_{C^{1+\beta}}\leq C||\phi||_{C^{1+\alpha}}||u-v||^\gamma
\end{align}

with $C=C(\alpha,\beta,||g||_{C^1},||\psi_u||_{C^1},||\psi_u||_{C^{1+\alpha}},L_0,L'_0,L_\alpha,L'_\alpha)$.
\bigskip

Now if $g$ also depends on $u\in\U$, computing $||(\W_u-\W_v)\phi||_{C^{1+\beta}}$ with $\phi\in C^{1+\alpha}$ would yield an additional term of the form $[(g(u,.)-g(v,.))\phi]\circ\psi(u,.)$, whose $C^{1+\beta}$ norm would be bounded by $C||\phi||_{C^{1+\alpha}}.||u-v||^\gamma$, with $C$ a constant. 
\\Thus, $u\in\mathcal{U}\mapsto\mathcal{L}_u\in L(C^{1+\alpha}(X),C^{1+\beta}(X))$ is (locally) $\gamma$-Hölder.
\bigskip

Let $\phi\in C^{1+\alpha}(W)$. The $C^{1}$ regularity of the inverse branches (w.r.t to $u$) allows one to consider the (partial) differential of $W$ with respect to $u$. Again, assume for the time being that $g$ does not depends on $u$.
Define $\chi_u: X\rightarrow L(\mathcal{B},TX)$ such that $D_u\psi_u=-\chi_u\circ\psi_u$, one gets :
\begin{align}
\partial_u\W(u,\phi)&= [Dg(\psi_u)\circ D_u\psi_u].\phi\circ\psi_u+g\circ\psi_u.[D\phi(\psi_u)\circ D_u\psi_u]
\end{align}

The previous formula defines a bounded operator $\partial_{u}\W\in L(\mathcal{B},L(C^{1+\alpha}(W),C^{\alpha}(W)))$, by virtue of lemma \ref{compisLipschitzlemma}.

One can easily extend the former to $\mathcal{L}_u$, and define a "partial differential" $\partial_u\mathcal{L}$, taking value in $L(\mathcal{B},L(C^{1+\alpha}(X),C^{\alpha}(X)))$. To what extend is it a "true" partial differential ? To answer that question one has to estimate $||\mathcal{L}(u_{0}+h,\phi)-\mathcal{L}(u_{0},\phi)-\partial_{u}\mathcal{.L}(u_{0},\phi).h||_{C^{\beta}}$, for $\phi\in C^{1+\alpha}(X)$

Let $x\in X$. 
One has
\begin{equation*}
[\W_{u_{0}+h}\phi-\W_{u_{0}}\phi-\partial_{u}\W(u_{0},\phi).h](x)=(I)+(II)+(III)
\end{equation*}
where
\begin{small}
\begin{align*}
(I)&=\phi(\psi(u_0,x))[g(\psi(u_0+h,x))-g(\psi(u_0,x))+Dg(\psi(u_0,x))\circ\chi_{u_{0}}(x).h]\\
(II)&=g(\psi(u_0,x))[\phi(\psi(u_0+h,x))-\phi(\psi(u_0,x))+D\phi(\psi(u_0,x))\circ\chi_{u_{0}}(x).h]\\
(III)&=[\phi(\psi(u_0+h,x))-\phi(\psi(u_0,x))][g(\psi(u_0+h,x))-g(\psi(u_0,x))]
\end{align*}
\end{small}

By lemma \ref{compisC1lemma},~\eqref{compisC1}, and lemma~\ref{compisLipschitzlemma},~\eqref{compisLipschitz} (I), (II) and (III) can be bounded as follows :
\begin{align*}
||(I)||_{C^\beta}&\leq C||\phi||_{C^{\beta}}||h||^{1+\gamma}||g||_{C^{1+\beta}}\\
||(II)||_{C^\beta}&\leq C||g||_{C^{\beta}}||h||^{1+\gamma}||\phi||_{C^{1+\beta}}\\
||(III)||_{C^\beta}&\leq C||h||^2.||\phi||_{C^{1+\beta}}||g||_{C^{1+\beta}}
\end{align*}

From the latter \footnote{From the previous bounds, one can even conclude that the map $u\in \mathcal{U}\longmapsto \mathcal{L}(u,\phi)\in C^{\beta}(X)$ is $C^{1+\gamma}$ for $\phi\in C^{1+\alpha}(X)$, which is precisely the conclusion drawn from the Taylor development at first order in Gouëzel-Liverani's paper (\cite[§8.1,~(8.3)]{GL06}).}, it is straightforward that
\begin{equation} \label{Taylor2}
 \mathcal{L}(u_{0}+h,\phi)-\mathcal{L}(u_{0},\phi)-\partial_{u}\mathcal{L}(u_{0},\phi).h=||h||_\mathcal{B}\epsilon(h,||g||_{C^{1+\beta}},||\phi||_{C^{1+\beta}})
\end{equation}
where $\epsilon(h,||g||_{C^{1+\beta}},||\phi||_{C^{1+\beta}})=\mathcal{O}(||h||_{\mathcal{B}}^\gamma)$.
\bigskip

Let us now show that $\mathcal{L}$ satisfies the Taylor expansion \eqref{NewTaylor} in the assumptions of theorem \ref{mainresult}. 
\\We start by recalling the following Taylor estimate, found in \cite{RDLL99}\footnote{We specifically refer to estimate (6.7) after theorem 6.10}:
 \\Letting E,F,G be Banach spaces, $\mathcal{U}\subset E$, $V\subset F$ be open sets, $0\leq\beta<\alpha<1$, and $(f,h)\in C^{1+\beta}(\mathcal{U},V)$ $(g,k)\in C^{1+\alpha}(V,G)$, one has
 \begin{equation}
 (g+k)\circ(f+h)=g\circ f+ k\circ f +[dg\circ f].h +R_{g,f}(h,k)
 \end{equation}
 where there exists some $0<\rho<1$ such that the remainder term $R_{g,f}(h,k)$ satisfies 
 \begin{equation}\label{compremainderestimate}
 ||R_{g,f}(h,k)||_{C^\beta}\leq C(||h||_{C^{1+\beta}}^{1+\rho}+||h||_{C^{1+\beta}}||k||_{C^{1+\alpha}})
 \end{equation}
 
This, together with the definition of $\partial_u\W_u$, yields for $(\phi,z)\in C^{1+\alpha}(W)$
 \begin{small}
	\begin{align}\label{Taylordvpmt}
	\notag \W_{u_0+h}(\phi+z)-\W_{u_0}(\phi)&-\partial_u\W(u_0,\phi).h-\W_{u_0}(z)\\
	&=D(g\phi)\circ\psi_{u_0}.(\psi_{u_0+h}-\psi_{u_0}-\partial_u\psi_{u_0}.h)+R_1(\psi_{u_0+h}-\psi_{u_0},g.z)
	\end{align}
\end{small}
where $R_1=R_{\phi,\psi_{u_0}}$ from \ref{compremainderestimate}. We start by bounding the first term. One has
 
\begin{equation}
 \psi_{u_0+h}-\psi_{u_0}-\partial_u\psi_{u_0}.h=\int_{0}^1[\partial_u\psi(u_0+th)-\partial_u\psi(u_0)].hdt  
\end{equation}
which leads us to estimate a term of the form $||df(\psi(u_0)).\int_{0}^1[\partial_u\psi(u_0+th)-\partial_u\psi(u_0)].hdt||_{C^\beta}$.
Following the trick used in the proof of lemma \ref{compisC1lemma}, we get
\begin{small}
	\begin{align}
 \notag	||df(\psi(u_0)).\int_{0}^1[\partial_u\psi(u_0+th)&-\partial_u\psi(u_0)].hdt||_{C^\beta}\\
 	&\leq [C_1||df||_{C^\beta}||\psi(u_0)||_{C^1}^\beta+C_2||df||_{\infty}]\frac{||h||^{1+\gamma}}{1+\gamma}
	\end{align}
\end{small}
 
 Now for $R_1$ we write, following estimate \eqref{compremainderestimate}:
 
 \begin{equation}
 ||R_1||_{C^\beta}\leq M[||h||^{1+\rho}+||h||.(C_1||z||_{C^{1+\alpha}}+C_0||z||_{C^\alpha})]
 \end{equation}
 
 with $C_1,~C_2$ depending on $\alpha,~||g||_{C^\alpha},~||g||_{C^{1+\alpha}}$.
 \bigskip
 
 Therefore, we obtained the following bound for \eqref{Taylordvpmt} :
 \begin{equation}
 M||h||^{1+\rho}+M'||h||^{1+\gamma}+C'_1||h||.||z||_{C^{1+\alpha}}+C'_2||h||.||z||_{C^\alpha}=[||h||+||z||_{C^\alpha}]\epsilon(h,z_{1+\alpha})
 \end{equation}
 where $z_{1+\alpha}$ is $z$ in $C^{1+\alpha}$ topology and $\epsilon(h,z_{1+\alpha})\underset{(h,z_{1+\alpha})\rightarrow 0}{\longrightarrow}0$ in $C^\beta(X)$.
 \bigskip
 
In the case of a weight $g$ depending on the parameter $u$, the partial derivative $\partial_u\W$ is given by 
\begin{equation}
\partial_u\W(u,\phi)= ([D_u(g)(u)]\phi)\circ\psi(u)+D_x(g\phi)\circ\psi(u).D_u\psi(u)
\end{equation}
Thus, the Taylor expansion at $(u_0,\phi)$ now has an additional term \begin{center}$[(g(u_0+h)-g(u_0)-D_u(g)(u_0).h)\phi]\circ\psi(u_0)$\end{center}
This term can be bounded (in $C^\beta$-norm), with upper bound of the form $C||g||_{C^{1+\alpha}}||h||^{1+\gamma}$, where $C=C(||\psi(u_0)||_{C^{1+\alpha}},||\phi||_{C^{1+\alpha}})$ is a constant, as outlined in lemma \ref{compisC1lemma}.
\\It follows that the transfer operator defined in~\eqref{transferopdef} also has a Taylor expansion of the form \eqref{NewTaylor}.

\begin{rem}\label{parameterremark}
The previous regularity results are given for $\mathcal{L}_u$ acting on the scale $(C^{1+\beta}(X),~C^\beta(X))$, $0<\beta<\alpha\leq1$. Following the method outlined in \cite{RDLL99}, and using theorem \ref{thetheorem}, one can show (by induction) that $\mathcal{L}_u$ acting on the scale $C^{k+\beta}(X),C^{k-j+\beta}(X)$ has a Taylor development of the form \eqref{NewTaylor2} at order j, with $0\leq j<k$ integers.
\end{rem}

\subsection{Hölder continuity of the spectral data : proof of theorem \ref{spectraldataisholder}}\label{holderspectraldatasection}

This section is devoted to establish theorem \ref{spectraldataisholder}, by  a direct argument. Note that this type of result is already known for a one-dimensional parameter, with previous works on spectral stability \cite{BY93,KL98}, or in the context of piecewise expanding maps of the interval \cite{K82}.
\medskip
 
Let $0\leq\beta<\alpha<1$, and $(T_{u})_{u\in\mathcal{\mathcal{U}}}$ be a family of $C^{1+\alpha}$ expanding maps, on a Riemann manifold X. Let $g:X\rightarrow\mathbb{R}$ be a positive~\footnote{Note that we only need the positivity of the weight to insure the simplicity of the maximal eigenvalue.} $C^{1+\alpha}$ function.
\\It follows from Ruelle theorem \cite{Ru89} that the transfer operator $(\mathcal{L}_{u})_{u\in \mathcal{U}}$ admits a spectral gap in $C^{1+\alpha}(X)$. Let $\lambda_{u}$ be the dominating eigenvalue of $\mathcal{L}_u$, $\phi_u\in C^{1+\alpha}(X)$ (resp $\ell_u\in(C^{1+\alpha}(X))^{'}$) be the right (resp left) eigenvector of $\mathcal{L}_u$ associated with $\lambda_u$, chosen so that $\langle \ell_u,\phi_u\rangle=1$. 
Let $F:\mathcal{U}\times C^{1+\alpha}(X)$, defined for $u\in \mathcal{U}$ and $\phi\not\in\ker\mathcal{L}_{u}^{*}\ell_{u_{0}}$, by
\begin{equation}\tag{\ref{normalisation}}
F(u,\phi)= \dfrac{\mathcal{L}_u\phi}{\langle\ell_{u_0},\mathcal{L}_u\phi\rangle}
\end{equation}
\bigskip

Note that F trivially inherits every regularity property of $(u,\phi)\in \mathcal{U}\times C_+^{1+\alpha}(X)^*\longmapsto\mathcal{L}_{u}\phi$, so in particular it is $\gamma$-Hölder in $u\in\U$ when considered as an operator from $C_+^{1+\alpha}(X)^*$ to $C_+^{1+\beta}(X)^*$. Hence the first point.
\medskip

The second item follows from the former remark and the fact that $\ell_{u_0}$ admits a bounded extension to $C^{1+\beta}(X)$, for every $0\leq\beta<\alpha$ (cf \cite{Ru89}). 
\\Let $\phi_{u}\in C_+^{1+\alpha}(X)^*$ be an eigenvector for $\lambda_{u}$, the dominating eigenvalue of $\mathcal{L}_{u}$. Then one has 

\begin{equation}
F(u,\phi_{u})=\dfrac{\lambda_{u}\phi_{u}}{\lambda_{u}\langle\ell_{u_{0}},\phi_{u}\rangle}=\dfrac{\phi_{u}}{\langle\ell_{u_{0}},\phi_{u}\rangle}
\end{equation}

For every $u\in \mathcal{U}$, fix a $\phi_{u}\in\ker(\lambda_{u}-\mathcal{L}_{u})$ such that $\langle\ell_{u_{0}},\phi_{u}\rangle=1$. Such a $\phi_{u}$ is unique in $\ker(\lambda_{u}-\mathcal{L}_{u})$ and verifies
\begin{equation}
F(u,\phi_{u})=\phi_{u}
\end{equation}

so that $F(u,.)$ has a unique fixed point $\phi_{u}$ in $C_+^{1+\alpha}(X)^*$ for every $u\in \mathcal{U}$.

Remark that for every $k\in\mathbb{N}^{*}$, for every $u\in \mathcal{U}$, every $\phi\not\in\ker((\mathcal{L}_{u}^{*})^{k}\ell_{u_{0}})$,
\begin{equation}\label{iterate}
F^{k}(u,\phi)=\dfrac{\mathcal{L}_{u}^{k}(\phi)}{\langle\ell_{u_{0}},\mathcal{L}_{u}^{k}(\phi)\rangle}
\end{equation}
by an immediate induction

Now note that, for every $k\in\mathbb{N}^{*}$, $u\in \mathcal{U}$,

\begin{align}
\phi(u)-\phi(u_{0})=F^{k}(u,\phi(u))-F^{k}(u_{0},\phi(u))+F^{k}(u_{0},\phi(u))-F^{k}(u_{0},\phi(u_{0}))\label{grossastuce3}
\end{align}

and that 
\begin{align}
F^{k}(u_{0},\phi(u))-F^{k}(u_{0},\phi(u_{0}))=\dfrac{\mathcal{L}_{u_0}^{k}(\phi(u))}{\langle\ell_{u_{0}},\mathcal{L}_{u_0}^{k}(\phi(u))\rangle}-\phi(u_{0})=\lambda_{u_{0}}^{-k}R_{u_{0}}^{k}(\phi(u)-\phi(u_{0}))
\end{align}

Recall that there is a $0<\sigma<1$ such that $||\lambda_{u_{0}}^{-k}R_{u_{0}}^{k}||_{C^{1+\beta}}\leq C\sigma^{k}$ (cf appendix \ref{expmapspectrum}), so that for $k$ large enough, one has

\begin{equation}
||F^{k}(u_{0},\phi(u))-F^{k}(u_{0},\phi(u_{0}))||_{C^{1+\beta}}\leq\dfrac{1}{2}||\phi(u)-\phi(u_{0})||_{C^{1+\beta}}
\end{equation}

From there,\eqref{grossastuce3} yields
\begin{align*}
||\phi(u)-\phi(u_{0})||_{C^{1+\beta}}&\leq C_{k,u}||u-u_{0}||^{\gamma}+ \dfrac{1}{2}||\phi(u)-\phi(u_{0})||_{C^{1+\beta}}\\
||\phi(u)-\phi(u_{0})||_{C^{1+\beta}}&\leq 2C_{k,u}||u-u_{0}||^{\gamma}
\end{align*}

where $C_{k,u}=||F^{k}(.,\phi(u))||_{C^{{1+\beta}}}$.
Thus, $u\in \mathcal{U}\longmapsto\phi(u)\in C^{1+\beta}(X)$ is $\gamma$-Hölder.

\subsection{Differentiability of the spectral data : proof of theorem \ref{linearresponse}}\label{diffspectraldatasection}
Let $0\leq\beta<\alpha<1$. This section is devoted to establish theorem \ref{linearresponse} by applying theorem \ref{mainresult} to the map $F$ from \eqref{normalisation} acting on the scale $(C^{1+\beta}(X),C^{\beta}(X))$.
\medskip

The first hypothesis, i.e existence, for every $u\in \mathcal{U}$, of a fixed point  $\phi_{u}$ for the map $F(u,.): C_+^{1+\alpha}(X)^*\rightarrow C_+^{1+\alpha}(X)^*$ from $\eqref{normalisation}$ and continuity of the map $u\in \mathcal{U}\mapsto\phi_{u}\in C^{1+\beta}(X)$, has already been addressed in theorem~\ref{spectraldataisholder}.
\bigskip

We now turn to assumption $(ii)$. We showed the perturbed Taylor development for $\mathcal{L}$ acting on $(C^{1+\beta}(X),C^{\beta}(X))$ in lemma~\ref{perturbedtransferop} : it immediately follows that $F$ acting on the scale $(C^{1+\beta}(X),C^{\beta}(X))$ satisfies the perturbed Taylor development \eqref{NewTaylor}.
\bigskip

We now check assumption $(iii)$. We start by remarking 
for every $z\in C^{1+\beta}(X)$,
\begin{equation}
Q_{u,\phi}.z=\frac{1}{\langle\ell_{u_{0}},\mathcal{L}(u,\phi)\rangle^{2}}[\mathcal{L}(u,z)\langle\ell_{u_{0}},\mathcal{L}(u,\phi)\rangle-\mathcal{L}(u,\phi)\langle\ell_{u_{0}},\mathcal{L}(u,z)\rangle]
\end{equation}
Thus, for $\phi=\phi_{u}$, we obtain

\begin{equation}
Q_{u,\phi_{u}}.z=\frac{1}{\lambda_{u}}(\mathcal{L}(u,z)-\langle\ell_{u_{0}},\mathcal{L}(u,z)\rangle\phi_{u})
\end{equation}

and for $u=u_{0}$ :

\begin{equation}
Q_{u_0,\phi_{u_0}}=\frac{1}{\lambda_{u_{0}}}\mathcal{L}(u_{0})-\Pi_{u_0}=\frac{1}{\lambda_{u_0}}R_{u_0}
\label{partial2}
\end{equation}

where $\Pi_{u_0}z=\langle\ell_{u_{0}},z\rangle\phi_{u_{0}}$, $z\in C^{1+\beta}(X)$ is the spectral projector on the (one-dimensional) eigenspace associated to $\lambda_{u_{0}}$. It is also noteworthy that the previous expression is independent of $\phi_{u_0}$.
\\From \eqref{partial2}, one sees that there is a $N\geq 1$ such that $||Q_{u_0}^N||_{C^\beta}\leq C\sigma^N$, for some $C>0$ and $\sigma\in(0,1)$ (cf appendix \ref{expmapspectrum},~\eqref{spectraldecomposition2}): therefore its Neumann series converges towards $(Id-Q_{u_0})^{-1}$.
\\This proves (iii) in the assumptions of theorem \ref{mainresult}.
\\We can therefore conclude that
\begin{center}
	If $\phi_u\in C^{1+\beta}(X)$, $u\in \mathcal{U}\longmapsto \phi_{u}\in C^{\beta}(X)$ is differentiable.
\end{center}

and that its differential satisfies
\begin{equation}
D_{u}\phi(u_{0})=(Id-Q_{u_0,\phi_{u_0}})^{-1}P_{u_0,\phi_{u_0}}
\end{equation}

Furthermore,
\begin{equation}
P_{u,\phi}=\dfrac{\partial_{u}\mathcal{L}(u,\phi)}{\langle\ell_{u_{0}},\mathcal{L}(u,\phi)\rangle}-\dfrac{\langle\ell_{u_{0}},\partial_{u}\mathcal{L}(u,\phi)\rangle}{\langle\ell_{u_{0}},\mathcal{L}(u,\phi)\rangle^{2}}\mathcal{L}(u,\phi)
\end{equation}

which simplifies, for $(u,\phi)=(u_{0},\phi_{u_{0}})$, to
\begin{align}\label{partial1}
P_{u_0,\phi_{u_0}}&=\frac{1}{\lambda_{u_{0}}}(\partial_{u}\mathcal{L}(u_{0},\phi_{u_{0}})-\langle\ell_{u_0},\partial_{u}\mathcal{L}(u_{0},\phi_{u_{0}})\rangle\phi_{u_{0}})\\
&=\frac{1}{\lambda_{u_{0}}}(Id-\Pi_{u_0})\circ\partial_{u}\mathcal{L}(u_{0},\phi_{u_{0}})
\end{align}
This, together with \eqref{partial2}, proves formula \eqref{linearresponseformula}.

\begin{corol}[Same setting as theorem \ref{linearresponse}]
The real valued map $u\in \mathcal{U}\longmapsto\lambda_{u}$ is differentiable
\end{corol}

\paragraph{Proof:} Let $u_0\in\mathcal{B}$, and $\U\subset\mathcal{B}$ be a neighborhood of $u_0$. Given the normalization chosen for $\ell_{u_0}$ and $\phi_u$ (cf §~\ref{mainresultstatement},~\eqref{normalisation}) for every $u\in\mathcal{U}$ one has
\begin{equation}
\lambda_u=\langle\ell_{u_0},\mathcal{L}(u,\phi_u)\rangle
\end{equation}
Thus, injecting \eqref{NewTaylor} and using the Hölder continuity (resp differentiability) of $u\in\U\mapsto\phi_u\in C^{1+\beta}(X)~(resp~C^{\beta}(X))$, one gets the desired conclusion.

\begin{corol}[Same setting as theorem \ref{linearresponse}]\label{Gmeasurediff}
Let $m_u$ be defined on $C^\beta(X)$ by $m_u(f)=\langle\ell_u,f\phi_u\rangle$. Then it is a Radon measure, and for every $f\in C^\beta(X)$, the map $u\in \mathcal{U}\longmapsto m_u(f)$ is $C^1$.
\end{corol}

\paragraph{Proof:}

By a standard positivity argument (see \cite{Ba00}) we extend continuously $\ell_u$ to $C^{0}(X)$. 
 It naturally follows that $m_u$ is a Radon measure. 
\bigskip

For $s\in D(0,1)\subset\mathbb{C},~u\in \mathcal{U}$ and $A\in C^{1+\alpha}(X)$, we introduce the parameter \\$\mathbf{u}=(s,u)\in D(0,1)\times \mathcal{U}\subset\mathbb{C}\times\mathcal{B}$ and the weighted transfer operator (with weight $e^g$, $g:X\rightarrow\mathbb{R}$) $\mathcal{L}_\mathbf{u}$ defined on $C^{1+\alpha}(X)$ by 
\begin{equation}
\mathcal{L}_\mathbf{u}\phi=\mathcal{L}_{s,u}=\mathcal{L}_u(e^{sA}\phi)
\end{equation}

Note that $\mathcal{L}_{s,u}$ is an analytical perturbation of $\mathcal{L}_u$ (at a fixed $u\in \mathcal{U}$). 
Hence, $\mathcal{L}_{s,u}$ also has a spectral gap for $s\in D(0,r)$, with $r=r(u)$ small enough (cf. \cite{Kat13}), and we will write $\lambda_{s,u},~\phi_{s,u}$ for its simple, maximal eigenvalue and the associated eigenvector (which is not necessarily a positive function, nor even a real valued one). 
\\It follows from Ruelle theorem \cite{Ru89} that $\lambda_{s,u}=e^{P(s,u)}$ with $P(s,u)$ the topological pressure associated with the dynamic $T_u$ and the weight $e^{sA+g}$. 
\\We now state a version of a well-known formula (cf. \cite{RuTF04}), connecting topological pressure and the expectation of the observable A under the Gibbs measure $m_u$, suited to our needs.

\begin{prop}
	Let $u\in \mathcal{U}$. The map $s\in D(0,r_u)\mapsto P(s,u)$ is analytical and one has
	\begin{equation}\label{pressureandmeasure}
	\partial_{s}P(0,u)=m_u(A)
	\end{equation}	
\end{prop}

\paragraph{Proof:} Fix $u\in \mathcal{U}$. For $s\in D(0,r)$, with $r=r(u)$ small enough, one can write \\$\mathcal{L}_{s,u}\phi_{s,u}=e^{P(s,u)}\phi_{s,u}$. The first statement follows from analytic perturbation theory, see \cite{Kat13}, as well as analyticity of $s\mapsto\ell_{s,u}$, with $\ell_{s,u}$ the eigenform for $\lambda_{s,u}$.

Furthermore, from the normalization $\langle\ell_{s,u},\phi_{s,u}\rangle=1$, one gets $\langle\ell_{s,u},\mathcal{L}_{s,u}\phi_{s,u}\rangle=e^{P(s,u)}$ and by differentiating this last equality with respect to s, one has 

\begin{equation}
\partial_sP(s,u)e^{P(s,u)}= \underset{(I)}{\underbrace{(\langle\partial_s\ell(s,u),\phi_{s,u}\rangle+\langle\ell_{s,u},\partial_s\phi_{s,u}\rangle)}e^{P(s,u)}}+\underset{(II)}{\underbrace{
\langle\ell_{s,u},\partial_{s}\mathcal{L}_{s,u}\phi_{s,u}\rangle}}
\end{equation}

From $\langle\ell_{s,u},\phi_{s,u}\rangle=1$, one gets $(I)=0$. 
\\Up to replace $A$ by $A\circ T$, $\partial_{s}\mathcal{L}_{s,u}\phi_{s,u}=A\mathcal{L}_{s,u}\phi_{s,u}=e^{P(s,u)}A\phi_{s,u}$, so that we get \\$(II)=e^{P(s,u)}\langle\ell_{s,u},A\phi_{s,u}\rangle$. Finally, one has, at $s=0$

\begin{equation}\tag{\ref{pressureandmeasure}}
\partial_sP(0,u)=\langle\ell_{u},A\phi_{u}\rangle=m_u(A)
\end{equation}

Fix a $u_0\in \mathcal{U}$: thus $\lambda_{0,u_0}=\lambda_{u_0}>0$.
\\ One easily has, for all $y\in X$,
\begin{center}
$\mathcal{L}_{s,u}\phi(y)=\underset{x\in T_u^{-1}y}{\sum}e^{sA(x)+g(x)}\phi(x)$
\end{center}
From theorem \ref{spectraldataisholder}, it holds that there is a neighborhood $D(0,r)\times B(u_0,\delta)$ such that $(s,u)\in D(0,r)\times B(u_0,\delta)$ implies $|\lambda_{s,u}-\lambda_{u_0}|\leq\frac{\lambda_{u_0}}{4}$. 
\\In particular, r is independent of u and $\lambda_{s,u}$ is a positive real number. Hence $P(s,u)$ is correctly defined, and continuous with respect to $u\in B(u_0,\delta)$, for $s\in D(0,r)$.
\bigskip

From theorem \ref{linearresponse}, it holds that there is a neighborhood $D(0,r')\times B(u_0,\delta')$ on which \\$(s,u)\longmapsto P(s,u)$ is $C^1$. In particular, $\partial_uP(s,u)$ exists and is continuous with respect to $u\in B(u_0,\delta')$ for $s\in D(0,r')$. Once again, $r'$ is \emph{a priori} independent of $u$.
\bigskip

From analytical perturbation theory, it holds that $s\in D(0,r")\longmapsto P(s,u)$ is analytical for $u\in B(u_0,\delta")$, where $r"=\min(r,r')$ and $\delta"=\min(\delta,\delta')$. Therefore, one can write, following Cauchy formula and \eqref{pressureandmeasure}
\begin{equation}
m_u(A)=\int_{\mathcal{C}(0,r")}\frac{P(s,u)}{s^2}ds
\end{equation}

where $\mathcal{C}(0,r")$ is the circle of radius $r"$ centered at 0. 
\\By Lebesgue's theorem, $u\in B(u_0,\delta")\longmapsto m_u(A)$ is a $C^1$ map. Up to a change in constants, this can be done for every $u_0\in U$, thus concluding this proof.

\appendix
\section{Spectrum of expanding maps on Hölder spaces}\label{expmapspectrum}

Recall that a $C^1$, expanding dynamic on a Riemann manifold X is a map $T:X\rightarrow X$ such that there exists a $\lambda>1$, and for every $x\in X$, every $v\in T_xX$, $||DT(x).v||\geq \lambda||v||$, where $TX$ is endowed with a norm field $(||.||_x)_{x\in X}$.
\\We recall a few useful properties of expanding maps in this setting:
\begin{prop}\label{propofexpmaps}
	Let $(X,g)$, $T$ be as above. Then
	\begin{itemize}
		\item[(i)] T is a local diffeomorphism at every $x\in X$.
		\item[(ii)] For every $y\in X$, $T^{-1}(\{y\})$ is a finite set.
		\item[(iii)] The set of $C^1$ expanding maps is open in the $C^1$-topology. Moreover, it is structurally stable.
	\end{itemize}
\end{prop}
The study of expanding maps started with the pioneering paper of Shub~\cite{Sh69}. One can find proof of the proposition claims in Shub's paper, or in the monograph \cite{KH97}. The study of their ergodic properties was started by \cite{KS69} where it is shown that every $C^2$ expanding map of a compact manifold has an invariant measure.
\medskip

Defining the (weighted) transfer operator associated to $(T,g)$ by 
\begin{equation}
\mathcal{L}\phi(x)=\sum_{y,Ty=x}g(y)\phi(y)
\end{equation}

where $g:X\rightarrow\mathbb{R}$, $C^{1}$ map, acting on the space $C^0(X)$, one can link statistical properties of the dynamic to spectral properties of $\mathcal{L}$ acting on an appropriate Banach space (\cite{Ba00,Liv02,Ba16}). As a result, the spectral picture of transfer operators for expanding maps has been thoroughly investigated, in the works of David Ruelle \cite{Ru89,Ru90},  Carlangelo Liverani\cite{Liv95,Liv02}, the 2000 monograph by Viviane Baladi \cite{Ba00}, or in a 2003 paper by Gundlash and Latushkin \cite{GuLa03}.

For example, SRB measures (which are physically relevant invariant measures, see \cite{LSY02}) and linear response formulas (first-order variation of the SRB measure w.r.t a real parameter) can be computed from spectral data of the transfer operator (\cite{Ru97,Ru97erratum,Liv02,Jiang12,Ba16}), decay of correlations can be linked to convergence of $\mathcal{L}^{n}$ towards its spectral projectors (\cite{Liv95,Ba00}).
\bigskip

The proper spectral setting is encapsulated in the notion \emph{spectral gap} : the operator $\mathcal{L}$ acting on the Banach space $\mathcal{B}$ has a spectral gap if :
\begin{itemize}
	\item There exists a simple, isolated eigenvalue $\lambda$ of maximal modulus, i.e $|\lambda|=\rho(\mathcal{L}_{|\mathcal{B}})$, called the \emph{dominating eigenvalue}.
	\item The rest of the spectrum is contained in a disk centered at 0 and of radius strictly smaller than $\rho(\mathcal{L}_{|\mathcal{B}})$.
\end{itemize}
In this case, one has the following decomposition :  

\begin{equation}\label{spectraldecomposition2}
\mathcal{L}\phi=\lambda\Pi(\phi)+R(\phi)
\end{equation}

In addition, the bounded operator R has the following property :
There exists $0<\sigma<1$, $C>0$ such that $||\lambda^{-n}R^{n}||_{\mathcal{B}}\leq C\sigma^{n}$.
\\Although $\mathcal{L}$ does not have nice spectral properties on $C^{0}(X)$ (\cite{Ru89}), a classical theorem of Ruelle (\cite{Ru89,Ru90}) shows that, assuming a little more regularity for the dynamic, the transfer operator admits a spectral gap on the Banach spaces $(C^{r}(X))_{r>0}$.

The proof relies on fine estimates on the (essential) spectral radius, established through \emph{Lasota-Yorke inequalities}. Those estimates were refined by Gundlash and Latushkin, in the paper \cite{GuLa03}, where they give an exact formula for the essential spectral radius of the transfer operator acting on $C^{r}(X)$ for $r\in\mathbb{R}_+$.
\bigskip 

A spectral gap can be obtained through other techniques, notably "cone contraction" based on abstract results of G.Birkhoff \cite{Bir57} : this approach was first applied in \cite{FS79} and successfully extended by C.Liverani \cite{Liv95}. Clear and complete account of those works can be found in the monographs by M.Viana \cite{Via97} or by V.Baladi\cite{Ba00}. Let us also mention the approach of Fan and Jiang \cite{FanJiang}.

\section{Estimates for compositions operators on Hölder spaces}\label{compestimates}
\paragraph{}
It is a well-established fact that $(C^{k+\alpha}(\Omega),||.||_{C^{k+\alpha}})$ is a Banach space. 
\\For $\Omega$ an open set in $\mathbb{R}^n$, and $0\leq\beta<\alpha<1$, one has the compact embedding :
\begin{center}
	$C^{k+\alpha}(\Omega)\Subset C^{k+\beta}(\Omega)$
\end{center}

The proof of this compact embedding relies on the Arzelà-Ascoli theorem and the following interpolation inequality :

\begin{theorem}
	Let E,F be Banach spaces, $\mathcal{U}\subset E$ an open subset. Let $0\leq\alpha<\beta<\gamma<1$ and $k\in\mathbb{N}$.
	\\Denote by $\mu=\frac{\gamma-\beta}{\gamma-\alpha}$. Then for every $f\in C^{k+\gamma}(\mathcal{U},F)$, one has
	\begin{equation}
	||f||_{C^{k+\beta}}\leq M_{\alpha}||f||_{C^{k+\alpha}}^\mu||f||_{C^{k+\gamma}}^{1-\mu}
	\end{equation}
\end{theorem}
We refer to \cite{RDLL99} for a proof.
\bigskip

The main object of this section is to address the regularity problem for composition operators: $g\longmapsto[f\longmapsto f\circ g]$ in Hölder spaces. An important inspiration for the results presented here is a paper by de la Llave and Obaya, \cite{RDLL99}, particularly the following result:

\begin{theorem}[\cite{RDLL99}, Prop 6.2, (iii)]\label{holderestimate}
	Let E,F,G be Banach spaces, and $\mathcal{U}\subset E$, $V\subset F$ open subsets. Let $k\geq 1$, $0\leq\gamma<1$ and $t=k+\gamma$. Let $s>t$ and $r\geq t$, and let $\mathcal{\mathcal{U}}\subset C^r(\mathcal{U},F)$. Then for every $g_1\in\mathcal{\mathcal{U}}$, there exists $\delta,\rho,M>0$, such that for every $f\in C^s(V,G)$, every $g_2\in C^r(\mathcal{U},F)$ which verifies $||g_1-g_2||_{C^r}\leq\delta$, one has $g_2\in\mathcal{\mathcal{U}}$, and
	
	\begin{equation}
	||f\circ g_1-f\circ g_2||_{C^t}\leq M||f||_{C^s}||g_1-g_2||_{C^r}^\rho
	\end{equation}
\end{theorem}
The estimates we establish in the following (lemmas~\ref{compisHölder2lemma},~\ref{compisLipschitzlemma},~\ref{compisC1lemma}~) are parameter variants of this theorem. They are used to prove lemma \ref{perturbedtransferop}, which in turn is key for using theorem \ref{mainresult} to prove theorem \ref{linearresponse}. 
\\In the first lemma, $g\longmapsto[f\longmapsto f\circ g]$ is Hölder continuous from $C^{1+\alpha}$ to $C^{1+\beta}$ with exponent $\gamma:=\alpha-\beta$

\begin{lemme}\label{compisHölder2lemma}
	Let $\mathcal{B},E,F,G$ be Banach spaces, $\mathcal{U}\subset\mathcal{B}$, $V\subset E$, $W\subset F$ be open domains. Let $0\leq\beta<\alpha<1$, $\psi\in C^0(\mathcal{U}\times V,W)$ such that for every $u\in \mathcal{U}$, $\psi_u=\psi(u,.)\in C^{1+\alpha}(V,W)$, and every $x\in V$, $u\longmapsto\psi(u,x)$ is Lipschitz continuous, $u\longmapsto D\psi_u(x)$ is $\alpha$-Hölder.
	Let $f\in C^{1+\alpha}(W,G)$.
	\\Denote by
	$$
	\left\{
	\begin{aligned}
	L_0=&\sup_{u\in \mathcal{U}}||\psi_u||_{Lip}~~L'_0=\sup_{x\in V}||\psi_.(x)||_{Lip}\\
	L_\alpha=&\sup_{u\in \mathcal{U}}||D\psi_u||_{C^\alpha}~~L'_\alpha=\sup_{x\in V}||D\psi_.(x)||_{\alpha}
	\end{aligned}
	\right.
	$$
	Let $u,v\in \mathcal{U}$. Then $f\circ\psi_u$, $f\circ\psi_v$ are $C^{1+\beta}$ maps, and we have
	\begin{equation}\label{compisHölder2}
	||f\circ\psi_u-f\circ\psi_v||_{C^{1+\beta}}\leq C||f||_{C^{1+\alpha}}||u-v||^\gamma
	\end{equation}
	with $C=C(\alpha,||f||_{C^1},||\psi_u||_{C^1},||\psi_u||_{C^{1+\alpha}},L_0,L'_0,L_\alpha,L'_\alpha)$
\end{lemme}

\paragraph{Proof:} We want to estimate $||f\circ\psi_u-f\circ\psi_v||_{C^{1+\beta}}=\max(||f\circ\psi_u-f\circ\psi_v||_{C^1},||D_x(f\circ\psi_u)-D_x(f\circ\psi_v)||_{C^\beta})$.
\\For $x\in V$, one has :
\begin{align*}
||Df(\psi_u(x))&\circ D\psi_u(x)-Df(\psi_v(x))\circ D\psi_v(x)||\\
&\leq ||Df(\psi_u(x))-Df(\psi_v(x))||.||D\psi_u(x)||+||Df(\psi_v(x))||.||D\psi_u(x)-D\psi_v(x)||\\
&\leq (||f||_{C^{1+\alpha}}||\psi_u||_{C^1}(L'_0)^\alpha+||f||_{C^1}L'_\alpha)||u-v||^\alpha.
\end{align*}

For the Hölder norm $||D_x(f\circ\psi_u)-D_x(f\circ\psi_v)||_{C^\beta}$, we have the following :
\\Let $x,x'\in V$, such that $||x-x'||\leq||u-v||$.
Then
\begin{align*}
&||Df(\psi_u(x))\circ D\psi_u(x)-Df(\psi_u(x'))\circ D\psi_u(x')||\\
&\leq||Df(\psi_u(x))-Df(\psi_u(x'))||.||D\psi_u(x)||+||Df(\psi_v(x))||.||D\psi_u(x)-D\psi_u(x')||\\
&\leq ||f||_{C^{1+\alpha}}||\psi_u(x)-\psi_u(x')||^\alpha+||f||_{C^1}L_\alpha||x-x'||^\alpha\\
&\leq (||f||_{C^{1+\alpha}}L_0^\alpha+||f||_{C^1}L_\alpha)||x-x||^\beta||u-v||^{\alpha-\beta}
\end{align*}

Similarly in the case $||x-x'||\geq||u-v||$, one has :
\begin{align*}
||Df(\psi_u(x))\circ D\psi_u(x)&-Df(\psi_v(x))\circ D\psi_v(x)||\\
&\leq (||f||_{C^1}L'_\alpha+||f||_{C^{1+\alpha}}||\psi_v||_{C^1}(L'_0)^\alpha)||u-v||^\gamma||x-x'||^\beta
\end{align*}
Thus, 
\begin{small}
	\begin{equation}
	|D_x(f\circ\psi_u)-D_x(f\circ\psi_v)|_{C^\beta}\leq 2(||f||_{C^{1+\alpha}}\max(L_0^\alpha,(L'_0)^\alpha)+||f||_{C^1}\max(L_\alpha,L'_\alpha))||u-v||^{\gamma}
	\end{equation}
\end{small}
and \eqref{compisHölder2} readily follows.
\paragraph{}
Note that the previous lemma yields Hölder continuity for $g\longmapsto[f\longmapsto f\circ g]$ from $C^{1+\alpha}$ to $C^{1+\beta}$, for $g\in C^{1+\alpha}(V)$.
One could easily follow the method outlined for the proof of theorem \ref{holderestimate} to establish our previous result from $C^{k+\alpha}(\Omega)$ to $C^{k+\beta}(\Omega)$, for every $k\geq 1$ and every $0\leq\beta<\alpha<1$.
\paragraph{}
One could ask what to expect for the composition operator from $C^{1+\alpha}$ to $C^\alpha$. We present the following estimate, a natural extension of the previous result
\begin{lemme}\label{compisLipschitzlemma}
	Let $\mathcal{B}$,E,F,G be Banach spaces, $\mathcal{U}\subset\mathcal{B}$, $V\subset E$, $W\subset F$ be open subsets. 
	\\Let $0\leq\alpha<1$ and $\psi\in C^{1+\alpha}(\mathcal{U}\times V,W)$, $f\in C^{1+\alpha}(W,G)$.
	\\Then for every $u_0\in \mathcal{U}$, and every $h\in\mathcal{B}$ such that $u_0+h\in \mathcal{U}$, the maps $f\circ\psi(u_0+h,.)$, $f\circ\psi(u_0,.)$ are $\alpha$-Hölder and one has the estimate:
	\begin{equation}\label{compisLipschitz}
	||f\circ\psi(u_0+h)-f\circ\psi(u_0)||_{C^\alpha}\leq C||f||_{C^{1+\alpha}}||h||_{\mathcal{B}}
	\end{equation} 
	with $C=C(\alpha,||\psi||_{C^1},||\psi||_{C^{1+\alpha}})$.	
\end{lemme}
\paragraph{Proof:} It is a straightforward consequence of the mean value theorem. Taking the $C^\alpha$-norm, one has for every $x\in V$.
\begin{equation}\label{meanvalueforLipestimate}
||f\circ\psi(u_0+h)-f\circ\psi(u_0)||_{C^\alpha}\leq||h||\int_{0}^1||Df(\psi(u_0+th))\circ D_u\psi(u_0+th)||_{C^\alpha}dt
\end{equation}
It is enough to establish the Lipschitz continuity that we wanted. Yet it is convenient to get a more precise estimate
of $||Df(\psi(u))\circ D_u\psi(u)||_{C^\alpha}$, for $u\in \mathcal{U}$.
\\Letting $x,x'\in W$, and taking the operator norm, one gets
\begin{align*}
||Df(\psi(u,x))&\circ D_u\psi(u,x)-Df(\psi(u,x'))\circ D_u\psi(u,x')||\\
&\leq ||Df(\psi(u,x))-Df(\psi(u,x'))||.||D_u\psi(u,x)||+||Df(\psi(u,x'))||.||D_u\psi(u,x)-D_u\psi(u,x')||\\
&\leq [||f||_{C^{1+\alpha}}||D_u\psi_u||_{C^0}||\psi_u||_{C^1}^\alpha+||f||_{C^1}||D_u\psi_u||_{C^\alpha}]||x-x'||^\alpha
\end{align*}
so that
\begin{small}
	\begin{equation}\label{holderestimate3}
	||Df(\psi(u))\circ D_u\psi(u)||_{C^\alpha}\leq ||f||_{C^{1+\alpha}}||D_u\psi_u||_{C^0}||\psi_u||_{C^1}^\alpha+||f||_{C^1}||D_u\psi_u||_{C^\alpha}
	\end{equation}
\end{small}
\bigskip

It is desirable to complete the previous lemmas with a differentiability result. In that spirit, we show the following
\begin{lemme}\label{compisC1lemma}
	Let $\mathcal{B}$,E,F,G be Banach spaces, $\mathcal{U}\subset\mathcal{B}$, $V\subset E$, $W\subset F$ be open subsets. 
	\\Let $0\leq\beta<\alpha<1$ and $\psi\in C^{1+\alpha}(\mathcal{U}\times V,W)$, $f\in C^{1+\alpha}(W,G)$.
	\\Denote by
	$$
	\left\{
	\begin{aligned}
	L_0=&\sup_{u\in \mathcal{U}}||\psi(u,.)||_{Lip}~~L'_0=\sup_{x\in\bar{\Omega}}||\psi(.,x)||_{Lip}\\
	L_{1,\alpha}=&\sup_{u\in \mathcal{U}}||D_u\psi(u,.)||_{C^\alpha}~~L'_{1,\alpha}=\sup_{x\in\Omega}||D_u\psi(.,x)||_{\alpha}
	\end{aligned}
	\right.
	$$
	Let $u_0\in \mathcal{U}$, and $h\in\mathbb{R}^d$ such that $u_0+h\in \mathcal{U}$. Then $f\circ\psi(u_0)$, $f\circ\psi(u_0+h)$, $D_u(f\circ\psi)(u_0)$ are $C^{\beta}$ maps, and we have
	\begin{equation}\label{compisC1}
	||f\circ\psi(u_0+h)-f\circ\psi(u_0)-D_u(f\circ\psi)(u_0).h||_{C^{\beta}}\leq C||f||_{C^{1+\alpha}}||h||^{1+\gamma}
	\end{equation}
	with $C=C(u_0,\alpha,||f||_{C^1},L_0,L'_0,L_{1,\alpha},L'_{1,\alpha})$
\end{lemme}

\paragraph{Proof:} Using the mean value theorem and taking the norm, one can write :
\begin{align}
\notag||f\circ\psi(u_0+h)&-f\circ\psi(u_0)-D_u(f\circ\psi)(u_0).h||_{C^\beta}\\
&\leq||h||\int_{0}^1||Df(\psi(u_0+th))\circ D_u\psi(u_0+th)-Df(\psi(u_0))\circ D_u\psi(u_0)||_{C^\beta}dt\label{meanvalueestimate}
\end{align}

To estimate $||Df(\psi(u_0+th))\circ D_u\psi(u_0+th)-Df(\psi(u_0))\circ D_u\psi(u_0)||_{C^\beta}$, we apply the same method we used to establish \eqref{compisHölder2}.
\\ Letting $x,x'\in V$, $u,v\in \mathcal{U}$, such that $||u-v||\leq||x-x'||$ one obtains :
\begin{small}
	\begin{align*}
	\frac{||Df(\psi(u,x))\circ D_u\psi(u,x)-Df(\psi(v,x))\circ D_u\psi(v,x)||}{||x-x'||^\beta}\leq(||f||_{C^{1+\alpha}}(L'_0)^\alpha||D_u\psi||_{\infty}+||f||_{C^1}L'_{1,\alpha})||u-v||^\gamma
	\end{align*}
\end{small}

Similarly, in the case $||x-x'||<||u-v||$
\begin{small}
	\begin{align*}
	\frac{||Df(\psi(u,x))\circ D_u\psi(u,x)-Df(\psi(u,x'))\circ D_u\psi(u,x')||}{||x-x'||^\beta}\leq[||f||_{C^{1+\alpha}}L_0^\alpha||D_u\psi||_{\infty}+||f||_{C^1}L_{1,\alpha}]||u-v||^\gamma
	\end{align*}
\end{small}

Finally, one has
\begin{small}
	\begin{equation*}
	\dfrac{||Df(\psi(u))\circ D_u\psi(u)-Df(\psi(v))\circ D_u\psi(v)||_{C^\beta}}{||u-v||^\gamma}\leq 2[||f||_{C^{1+\alpha}}||D_u\psi||_{\infty}\max(L_0,L'_0)^\alpha+||f||_{C^1}\max(L_{1,\alpha},L_{1,\alpha}')]
	\end{equation*}
\end{small}

Injecting this last estimate in \eqref{meanvalueestimate}, one gets the following :
\begin{align*}
||f\circ\psi(u_0+h)-f\circ\psi(u_0)-D_u(f\circ\psi)(u_0).h||_{C^\beta}&\leq||h||\int_{0}^1C||f||_{C^{1+\alpha}}t^\gamma||h||^\gamma dt\\
&=C||h||^{1+\gamma}||f||_{C^{1+\alpha}}\frac{1}{1+\gamma}
\end{align*}
which gives the promised result with $C'=\frac{C}{1+\gamma}$.

\section{An elementary example}\label{elementary}
Let $I=[-1,1]$, and consider the Banach space $C^{0}(I)$. Let $0<\epsilon<1$, and define the family of maps $(F_{u})_{u\in [-\epsilon,\epsilon]}$ by :
\begin{equation}
F_{u}(\phi)(t)=\frac{1}{2}\phi(\frac{t+u}{2})+g(t,u)
\end{equation}
with $g:I\times [-\epsilon,\epsilon]\rightarrow \mathbb{R}$ a non-zero $C^{1}$ map, such that $g(t,.)\in B_{C^{\alpha}}(0,\dfrac{1}{2})$. 
Being a contraction of $C^{0}(I)$, $F_{u}$ admits a fixed point, say $\phi_{u}$, by the Banach contraction principle. But what about the regularity of $u\in I\mapsto\phi_{u}\in C^{0}(I)$ ?
On the $C^{0}$ space, $u\mapsto F_{u}$ is not even continuous. Nevertheless, if we consider the same operator $F(u,.):C^{\alpha}(I)\rightarrow C^{\alpha}(I)$, with $\alpha\in(0,1)$, and the immersion $\mathcal{I}:C^{\alpha}(I)\rightarrow C^{0}(I)$, we see that :

\begin{align*}
|(F(u,\phi)-F(u,\psi))(t)-(F(u,\phi)-F(u,\psi))(t')|\leq\frac{1}{2^{1+\alpha}}||\phi-\psi||_{C^{\alpha}}|t-t'|^{\alpha}
\end{align*}

It follows that,

\begin{equation*}
||F(u,\phi)-F(u,\psi)||_{C^{\alpha}}\leq \frac{1}{2^{1+\alpha}}||\phi-\psi||_{C^{\alpha}}
\end{equation*}

and $F(u,.):C^{\alpha}(I)\rightarrow C^{\alpha}(I)$ is a contraction; by the Banach contraction principle, this map has a fixed point in $C^{\alpha}(I)$ for all $u\in I$, say $\phi(u)$. Note also that 
\begin{align*}
|F(u,\phi)(t)-F(u',\phi)(t)|\leq(\frac{1}{2^{1+\alpha}}||\phi||_{C^{\alpha}}+||g||_{C^{\alpha}})|u-u'|^{\alpha}
\end{align*}
so $u\in I\mapsto F(u,.)\in C^{0}(I)$ is a $\alpha$-Hölder map. Finally, 

\begin{align*}
||\phi(u)-\phi(u')||_{C^{0}}&=||F(u,\phi(u))-F(u',\phi(u'))||_{C^{0}}\\
&=||F(u,\phi(u))-F(u',\phi(u))+F(u',\phi(u))-F(u',\phi(u'))||_{C^{0}}\\
&\leq(\frac{1}{2^{1+\alpha}}||\phi(u)||_{C^{\alpha}}+||g||_{C^{\alpha}})|u-u'|^{\alpha}+ \frac{1}{2}||\phi(u)-\phi(u')||_{C^{0}}
\end{align*} 

Hence

\begin{equation*}
||\phi(u)-\phi(u')||_{C^{0}}\leq (\frac{1}{2^{\alpha}}||\phi(u)||_{C^{\alpha}}+2||g||_{C^{\alpha}})|u-u'|^{\alpha} 
\end{equation*}

and $u\mapsto \phi(u)\in \mathcal{I}(C^{\alpha})$ is locally $\alpha$-Hölder.

\begin{footnotesize}
\bibliographystyle{plain}
\bibliography{biblio}

\begin{thebibliography}{10}

\bibitem{Ba00}
Viviane Baladi.
\newblock {\em Positive Transfer Operator and Decay of Correlations}.
\newblock World Scientific, 2000.

\bibitem{Ba14}
Viviane Baladi.
\newblock Linear response, or else.
\newblock {\em arXiv preprint arXiv:1408.2937}, 2014.

\bibitem{Ba16}
Viviane Baladi.
\newblock {\em Dynamical zeta functions and dynamical determinants for
  hyperbolic maps. A functional approach.}
\newblock A paraître, Springer, 2016.

\bibitem{BaladiSmania08}
Viviane Baladi and Daniel Smania.
\newblock Linear response formula for piecewise expanding unimodal maps.
\newblock {\em Nonlinearity}, 21(4):677, 2008.

\bibitem{BaladiSmaniaalt10}
Viviane Baladi and Daniel Smania.
\newblock Alternative proofs of linear response for piecewise expanding
  unimodal maps.
\newblock {\em Ergodic Theory and Dynamical Systems}, 30(1):1--20, 2010.

\bibitem{BaladiSmania12}
Viviane Baladi and Daniel Smania.
\newblock Linear response for smooth deformations of generic nonuniformly
  hyperbolic unimodal maps.
\newblock {\em Ann. Sci. {\'E}c. Norm. Sup{\'e}r}, 45(6):861--926, 2012.

\bibitem{BaladiTodd16}
Viviane Baladi and Mike Todd.
\newblock Linear response for intermittent maps.
\newblock {\em Communications in Mathematical Physics}, 347(3):857--874, 2016.

\bibitem{BY93}
Viviane Baladi and L-S Young.
\newblock On the spectra of randomly perturbed expanding maps.
\newblock {\em Communications in Mathematical Physics}, 156(2):355--385, 1993.

\bibitem{Bir57}
Garrett Birkhoff.
\newblock Extensions of {J}entzsch's theorem.
\newblock {\em Transactions of the American Mathematical Society},
  85(1):219--227, 1957.

\bibitem{RDLL99}
Rafael De~la Llave and Rafael Obaya.
\newblock Regularity of the composition operator in spaces of {H}{\"o}lder
  functions.
\newblock {\em Discrete and Continuous Dynamical Systems}, 5:157--184, 1999.

\bibitem{Dol04}
Dmitry Dolgopyat.
\newblock On differentiability of {SRB} states for partially hyperbolic
  systems.
\newblock {\em Inventiones mathematicae}, 155(2):389--449, 2004.

\bibitem{FanJiang}
Aihua Fan and Yunping Jiang.
\newblock Spectral theory of transfer operators.
\newblock {\em Jiang Y, Wang Y. Complex Dynamics and Related Topics. New
  Studies in Advanced Mathematics}, 5:63--128, 2004.

\bibitem{FS79}
P~Ferrero and B~Schmitt.
\newblock Ruelle’s {P}erron-{F}robenius theorem and projective metrics.
\newblock {\em Coll.Math.Soc}, 27:333--336, 1979.

\bibitem{GL08}
S{\'e}bastien Gou{\"e}zel, Carlangelo Liverani, et~al.
\newblock Compact locally maximal hyperbolic sets for smooth maps: fine
  statistical properties.
\newblock {\em Journal of Differential Geometry}, 79(3):433--477, 2008.

\bibitem{G10}
Sébastien Gouëzel.
\newblock Characterization of weak convergence of {B}irkhoff sums for
  {G}ibbs-{M}arkov maps.
\newblock {\em Israël Journal of Mathematics}, 180:1--41, December 2010.

\bibitem{GL06}
Sébastien Gouëzel and Carlangelo Liverani.
\newblock Banach spaces adapted to {A}nosov systems.
\newblock {\em Ergodic Theory and Dynamical systems}, 26:189--217, February
  2006.

\bibitem{GuLa03}
VM~Gundlach and Yu~Latushkin.
\newblock A sharp formula for the essential spectral radius of the {R}uelle
  transfer operator on smooth and {H\"o}lder spaces.
\newblock {\em Ergodic Theory and Dynamical Systems}, 23(1):175--191, 2003.

\bibitem{HairerMajda10}
Martin Hairer and Andrew~J Majda.
\newblock A simple framework to justify linear response theory.
\newblock {\em Nonlinearity}, 23(4):909, 2010.

\bibitem{hamilton1982inverse}
Richard~S Hamilton.
\newblock The inverse function theorem of {N}ash and {M}oser.
\newblock {\em Bulletin of the American Mathematical Society}, 7(1):65--222,
  1982.

\bibitem{Jiang12}
Miaohua Jiang.
\newblock Differentiating potential functions of {S}{R}{B} measures on
  hyperbolic attractors.
\newblock {\em Ergodic Theory and Dynamical Systems}, 32(4):1350--1369, 2012.

\bibitem{Kat13}
Tosio Kato.
\newblock {\em Perturbation theory for linear operators}, volume 132.
\newblock Springer Science \& Business Media, 2013.

\bibitem{KH97}
Anatole Katok and Boris Hasselblatt.
\newblock {\em Introduction to the modern theory of dynamical systems},
  volume~54.
\newblock Cambridge university press, 1997.

\bibitem{K82}
Gerhard Keller.
\newblock Stochastic stability in some chaotic dynamical systems.
\newblock {\em Monatshefte f{\"u}r Mathematik}, 94(4):313--333, 1982.

\bibitem{KS69}
K~Krzy{\.z}ewski and W\_ Szlenk.
\newblock On invariant measures for expanding differentiable mappings.
\newblock In {\em The Theory of Chaotic Attractors}, pages 37--46. Springer,
  1969.

\bibitem{Liv95}
Carlangelo Liverani.
\newblock Decay of correlations.
\newblock {\em Annals of Mathematics}, 142(2):239--301, September 1995.

\bibitem{Liv02}
Carlangelo Liverani.
\newblock Invariant measures and their properties : a functional analytic point
  of view.
\newblock In {\em Dynamical Systems.Part II: Topological Geometrical and
  Ergodic Properties of Dynamics}. Scuola Norm. Sup., Pisa, Pubbl. Cent. Ric.
  Mat. Ennio Giorgi, 2003.

\bibitem{KL98}
Carlangelo Liverani and Gerhard Keller.
\newblock Stability of the spectrum for transfer operators.
\newblock {\em Ann. Scuola Norm. Sup. Pisa Cl. Sci.(4) vol}, 28:141--152, 1998.

\bibitem{Ru89}
David Ruelle.
\newblock Thermodynamic formalism for expanding maps.
\newblock {\em Comm. Math. Phys.}, 1989.

\bibitem{Ru90}
David Ruelle.
\newblock An extension of the theory of fredholm determinants.
\newblock {\em Inst. Hautes Etudes Sci. Publ. Math}, 72:175--193, 1990.

\bibitem{Ru97}
David Ruelle.
\newblock Differentiation of {SRB} states.
\newblock {\em Communications in Mathematical Physics}, 187:227--241, July
  1997.

\bibitem{Ru98}
David Ruelle.
\newblock Nonequilibrium statistical mechanics near equilibrium: computing
  higher-order terms.
\newblock {\em Nonlinearity}, 11(1):5, 1998.

\bibitem{Ru97erratum}
David Ruelle.
\newblock Differentiation of {S}{R}{B} states: correction and complements.
\newblock {\em Communications in mathematical physics}, 234(1):185--190, 2003.

\bibitem{RuTF04}
David Ruelle.
\newblock {\em Thermodynamic formalism: the mathematical structure of
  equilibrium statistical mechanics}.
\newblock Cambridge University Press, 2004.

\bibitem{Rugh08}
Hans~Henrik Rugh.
\newblock On the dimensions of conformal repellers. {R}andomness and
  {P}arameter dependency.
\newblock {\em Annals of Mathematics}, pages 695--748, 2008.

\bibitem{Rugh10}
Hans~Henrik Rugh.
\newblock Cones and gauges in complex spaces: Spectral gaps and complex
  {P}erron-{F}robenius theory.
\newblock {\em Annals of mathematics}, pages 1707--1752, 2010.

\bibitem{Sh69}
M~Shub.
\newblock Endomorphisms of compact differentiable manifolds.
\newblock {\em Amer. J. Maths.}, 91:175--199, 1969.

\bibitem{Via97}
Marcelo Viana.
\newblock {\em Stochastic dynamics of deterministic systems}, volume~21.
\newblock IMPA Rio de Janeiro, 1997.

\bibitem{LSY02}
Lai-Sang Young.
\newblock What are {S}{R}{B} measures, and which dynamical systems have them?
\newblock {\em Journal of Statistical Physics}, 108(5):733--754, 2002.

\end{thebibliography}
\end{footnotesize}

\end{document}